\documentclass[10pt]{article}

\usepackage{stmaryrd}
\usepackage{mathtools}
\usepackage[margin=1in]{geometry}
\usepackage{xfrac}
\usepackage{epsfig}
\usepackage{amssymb,amsmath}
\usepackage{enumerate}
\usepackage{amsfonts}
\usepackage{graphics}
\usepackage{multirow}
\usepackage{mathrsfs}
\usepackage{comment}
\epsfverbosetrue
\usepackage {amssymb,rotating,graphicx,cite, calc,color,dsfont,epsfig,slashbox,relsize}
\usepackage{epsfig}
\usepackage{graphicx}
\usepackage{color}
\usepackage{xcolor}
\usepackage{stfloats}
\usepackage{mathrsfs}
\usepackage{todonotes}
\usepackage{cite}
\usepackage[ruled,vlined]{algorithm2e}
\def\zd{\mathcal{Z}_d}
\definecolor{blue}{RGB}{0,0,230}
\usetikzlibrary{arrows,calc,decorations.markings,math,arrows.meta}
\newcommand{\pf}{{ \hspace{.1in}{  {\it Proof:}} \ }}

\newcommand{\tinf}{t\rightarrow\infty}

\newcommand{\Rn}{\mathbb{R}^n}
\newcommand{\Rd}{\mathbb{R}^m}
\newcommand{\Rdn}{(\mathbb{R}^{m})^n}
\def\R{\mathbb{R}}

\makeatletter
\newcommand*{\rom}[1]{\expandafter\@slowromancap\romannumeral #1@}
\makeatother

\newcommand{\E}{\mathbb{E}}

\newcommand{\N}{\mathbb{N}}
\newcommand{\F}{\mathcal{F}}
\newcommand{\EE}{\mathcal{E}}
\newcommand{\ED}{\mathcal{D}}
\newcommand{\EH}{\mathcal{H}}
\newcommand{\G}{\mathcal{G}}

\newcommand{\diam}{{\rm d}}
\newcommand{\diamm}{{\rm diam}}
\newcommand{\bx}{{\bf x}}
\newcommand{\by}{{\bf y}}
\newcommand{\bh}{{y}}
\newcommand{\bg}{{\bf g}}
\newcommand{\bv}{{ v}}

\def\u{{\bf u}}

\newcommand{\qed}{\hfill \rule{1.1ex}{1.1ex} \\}
\newtheorem{theorem}{Theorem}

\newtheorem{lemma}{Lemma}

\newtheorem{example}{Example}
\newtheorem{assumption}{Assumption}

\newcommand{\lrangle}[1]{\left\langle #1 \right\rangle}

\title{
Distributed Optimization Over Dependent Random Networks}
\author{Adel Aghajan and Behrouz Touri\thanks{This work is supported by the NSF Grant 1913131. }\\ University of California, San Diego, email: \textit{\{adaghaja,btouri\}@ucsd.edu}}
\begin{document}

\date{}
\maketitle
\sloppy
% \begin{abstract}

% \end{abstract}
\begin{abstract}
We study the averaging-based distributed optimization solvers over random networks. We show a general result on the convergence of such schemes using weight-matrices that are row-stochastic almost surely and column-stochastic in expectation for a broad class of dependent weight-matrix sequences. In addition to implying many of the previously known results on this domain, our work shows the robustness of distributed optimization results to link-failure. Also, it provides a new tool for synthesizing distributed optimization algorithms. {To prove our main theorem, we establish new results on the rate of convergence analysis of averaging dynamics over (dependent) random networks. These secondary results, along with the required martingale-type results to establish them, might be of interest to a broader research endeavors in distributed computation over random networks.  }
\end{abstract}
\section{Introduction}
 
Distributed optimization has received increasing attention in recent years due to its applications in
 distributed control of robotic networks \cite{bullo2009distributed}, study of opinion dynamics in
social networks \cite{hegselmann2002opinion},  distributed estimation and signal processing \cite{rabbat2004distributed,cattivelli2010diffusion,stankovic2011decentralized}, and power networks \cite{dominguez2011distributed,dominguez2012decentralized,cherukuri2015distributed}.

In distributed optimization, we are often interested in finding an optimizer of a decomposable function $F({z})=\sum_{i=1}^{n}f_{i}({z})$ such that  $f_i(\cdot)$s are distributed through a network of $n$ agents, i.e., agent $i$ only knows $f_i(\cdot)$, and we are seeking to solve this problem without sharing the local objective function $f_i(z)$. Therefore, the goal is to find distributed dynamics over time-varying networks that, asymptotically, all the nodes agree on an optimizer of $F(\cdot)$. 
%Therefore, we use consensus algorithms that make the value of every node converge to a similar value. To do so, we need row-stochastic non-negative matrices, i.e., the summation of each row of the matrix is equal to one. 

The most well-know algorithm that achieves this is, what we refer to as, the averaging-based distributed optimization solver where each node maintains an estimate of an optimal point, and at each time step, each node computes the average of the estimates of its neighbors and performs (sub-)gradient descent on its local objective function \cite{nedic2009distributed}. However, in order for such an algorithm to converge, the corresponding weight matrices should be doubly stochastic.  While making a row-stochastic or a column-stochastic matrix is easy, this is not the case for doubly stochastic matrices. Therefore, it is often assumed that such a doubly stochastic matrix sequence is given.

A solution to overcome this challenge is to establish more complicated distributed algorithms that effectively \textit{reconstruct} the average-state distributively. The first algorithm in this category was proposed in \cite{rabbat}, which is called subgradient-push (or push-sum), and later was extended for time-varying networks \cite{nedic2014distributed}. In this scheme, the weight matrices are assumed to be column-stochastic, and through the use of auxiliary state variables the approximate average state is reconstructed. 
% Unlike the row-stochastic weight matrices that are made at the receivers, column-stochastic matrices are pushed by senders to receivers. Thus, if one link fails, the respective weight matrix is not a column-stochastic matrix anymore. Therefore, such algorithms are prone to link-failures. 
Another scheme in this category that works with row-stochastic matrices, but does not need the column-stochastic assumption, is proposed in \cite{mai2016distributed,xi2018linear}. However, to use  this scheme, every node needs to be assigned and know its unique label. Assigning those labels distributively is also another challenge in this respect. In addition, both these schemes invoke  division operation which results in theoretical challenges in establishing their stability in random networks \cite{rezaeinia2019push,rezaeinia2020distributed}. 

{Another solution to address this challenge is to use gossip \cite{7165636} and broadcast gossip \cite{5585721} algorithms over random networks. The weight matrices of gossip algorithms are row-stochastic and in-expectation column-stochastic. This fact was generalized in \cite{7849184}, where it is proven that it is sufficient to have row-stochastic weight matrices, that are column-stochastic in-expectation. In all the above works on distributed optimization over random networks, all weight matrices are assumed to be independent and identically distributed (i.i.d.). In \cite{lobel2011distributed}, the broader class of random networks, which is Markovian networks was studied in distributed optimization; however,  weight matrices were assumed to be doubly stochastic almost surely.   
Also, our work is closely related to the existing works on distributed averaging on random networks \cite{consensus1,consensus2,consensus3,touri2014endogenous}.
% , and distributed optimization \cite{lobel2010distributed,nedic2007rate,nedic2009distributed,do1,do2,do3,nedic2014distributed,mai2016distributed}.

In this paper, we study distributed optimization over random networks, where the randomness is not only time-varying but also, possibly, dependent on the past. Under the standard assumptions on the local objective functions and step-size sequences for the gradient descent algorithm, we show that the averaging-based distributed optimization solver at each node converges to a global optimizer almost surely if the weight matrices are row-stochastic almost surely, column-stochastic {in-expectation}, and satisfy certain connectivity assumptions.   }
% To the best of our knowledge this work is the first to recognize this feature of the averaging-based distributed optimization solvers, and we believe that it has many implications on other extensions of this algorithm (a.e.\ non-convex settings \cite{nc1,nc2}, and accelerated algorithms \cite{accelarated0,accelerated1,accelarated2,accelarated3}). 

% The extensive list of the previous works on distributed optimization and distributed averaging is well beyond the scope of this paper. However, the main related work to our study is the earlier work \cite{lobel2010distributed}, where the weight matrices are assumed to be random, but still, the doubly stochastic assumption on the weight matrices is assumed almost surely.  

% The main contribution of the current work is as follows: we show that under the standard assumptions on the local objective functions, step-size sequences for the gradient descent algorithm, and the connectivity of the networks, the averaging-based distributed optimization solver at each node converges to a global optimizer almost surely if the sequence of row-stochastic weight matrices are independent and are doubly stochastic \textit{in-expectation}. To the best of our knowledge this work is the first to recognize this feature of the averaging-based distributed optimization solvers, and we believe that it has many implications on other extensions of this algorithm (a.e.\ non-convex settings \cite{nc1,nc2}, and accelerated algorithms \cite{accelarated0,accelerated1,accelarated2,accelarated3}). 

The paper is organized as follows: we conclude this section by introducing  mathematical notations that will be used subsequently. In Section \ref{sec:ProblemStatement}, we formulate the problem of interest and state the main result of this work and discuss some of its immediate consequences. To prove the main result, first we study the behavior of the distributed averaging dynamics over random networks in Section~\ref{sec:ConsensusWithoutPerturbation}. Then, in Section \ref{sec:ConsensusWithPerturbation}, we extent this analysis to the dynamics with arbitrary control inputs. Finally, the main result is proved in Section \ref{sec:Convergence}. We conclude this work in Section~\ref{Conclusion}.

\textbf{Notation and Basic Terminology:} The  following  notation  will be  used  throughout the paper. We let $[n]\triangleq\{1,\ldots,n\}$. We denote the space of real numbers by $\R$  and natural (positive integer) numbers  by $\mathbb{N}$. We denote the space of $n$-dimensional real-valued vectors by $\R^n$. In this paper, all vectors are assumed to be column vectors. The transpose of a vector $x\in \R^n$ is denoted by $x^T$.
For a vector $x\in \R^n$, $x_i$ represents the $i$th coordinate of $x$. %, except where for notational convenience, we denote $e_1^n,\ldots,e_n^n$ as the standard basis vectors of $\Rn$. 
We denote the all-one vector in $\Rn$ by $e^n=[1,1,\ldots,1]^T$. We drop the superscript $n$ in $e^n$ whenever the dimension of the space is understandable from the context. %For two $m\times n$ matrices $A$ and $B$, we use $A\geq B$ when $a_{ij}\geq b_{ij}$ for all $i\in [m]$ and all $j\in [n]$.
  A non-negative matrix $A$ is a row-stochastic (column-stochastic) matrix if $Ae=e$ ($e^TA= e^T$).
  
Let $(\Omega,\mathcal{F},\text{Pr})$ be a
probability space and let $\{W(t)\}$ be a chain of random
 matrices, i.e., for all $t\geq 0$ and $i,j\in[n]$, $w_{ij}(t):\Omega\to\R$ is a Borel-measurable function.
 For random vectors (variables) { $\bx(1),\ldots,\bx(t)$, we denote the sigma-algebra generated by these random variables by $\sigma(\bx(0),\ldots,\bx(t))$.} We say that $\{\F(t)\}$ is a filtration for $(\Omega,\mathcal{F})$ if $\F(0)\subseteq \F(1)\subseteq \cdots \subseteq \F$. Further, we say that a random process $\{V(t)\}$ (of random variables, vectors, or matrices) is adapted to $\{\F(t)\}$ if $V(t)$ is measurable with respect to $\F(t)$.

Throughout this paper we mainly deal with directed graphs. A directed graph $\G=([n],\EE)$ (on $n$ vertices) is defined by a vertex set (identified by) $[n]$ and an edge set $\EE\subset [n]\times [n]$.
{A graph $\G=([n],\EE)$ has a spanning directed rooted tree if it has a vertex $r\in [n]$ as \textit{a} root such that there exists a (directed) path from $r$ to every other vertex in the graph.}
For a matrix $A=[a_{ij}]_{n\times n}$, the associated directed graph with parameter $\gamma>0$ is the graph $\G^\gamma(A)=([n],\EE^\gamma(A))$ with the edge set $\EE^\gamma(A)=\{(j,i)\mid i,j\in [n],a_{ij}>\gamma \}$. Later, we fix the value $0<\gamma<1$ throughout the paper and hence, unless otherwise stated, for notational convenience, we use $\G(A)$ and $\EE(A)$ instead of $\G^\gamma(A)$ and $\EE^\gamma(A)$. 

{ The function $f:\Rd\to\R$ is convex if for all $x,y\in\Rd$ and all $\theta\in[0,1]$, 
\begin{align*}
    f(\theta x+(1-\theta)y)\leq \theta f(x)+(1-\theta)f(y).
\end{align*}
We say that $g\in\Rd$ is a subgradient of the function $f(\cdot)$ at $\hat{x}$ if  for all $x\in\Rd$,
% \begin{align*}
 $f(x)-f(\hat{x})\geq \lrangle{g,x-\hat{x}},$  
% \end{align*}
where $\lrangle{u_1,u_2}=u_1^Tu_2$ is the standard inner product in $\Rd$.}
The set of all subgradients of $f(\cdot)$ at $x$ is denoted by $\nabla f(x)$.
For a convex function $f(\cdot)$, $\nabla f(x)$ is not empty for all $x\in\Rd$ (see e.g., Theorem 3.1.15 in \cite{nesterov2018lectures}).  Finally, for convenience and due to the frequent use of $\ell_\infty$ norm in the study of averaging dynamics, we use $\|\cdot\|$ to denote the $\ell_\infty$ norm $\|x\|\triangleq \max_{i\in [m]}|x_i|$. 
\section{Problem Formulation and Main Result}\label{sec:ProblemStatement}
In this section, we discuss the main problem and the main result  of  this  work.  The  proof  of  the  result  is provided in the subsequent sections.
\subsection{General Framework}
Consider a communication network with $n$ nodes or agents such that  node $i$ has the cost function $f_i:\Rd\to\R$. Let $F({z})\triangleq\sum_{i=1}^{n}f_{i}({z})$. The goal of this paper is to solve   
\begin{align}\label{eqn:mainProblem}
    \arg\min_{{z}\in\Rd}F({z})
\end{align}
distributively with the following assumption on the objective function.
\begin{assumption}[Assumption on the Objective Function]\label{assump:AssumpOnFunction}
We assume that: 
\begin{enumerate}[(a)]
    \item All $f_i({z})$ are convex functions over $\Rd$. 
    \item The optimizer set $\mathcal{Z}\triangleq\arg\min_{{z}\in\Rd}F({z})$ is non-empty.
    \item The subgradients $f_i({z})$s are uniformly upper bounded, i.e., for all $g\in \nabla f_i(z)$, $\|g\|\leq L_i$ for all $z\in \Rd$ and all $i\in[n]$. {We let $L\triangleq\sum_{i=1}^nL_i$.}
\end{enumerate}
\end{assumption}
%In this paper, for matrices, we use capital letters, such as $W,A$, except matrices $\bx(t),\bg(t)\in\mathbb{R}^{n\times d}$. Also, we use $\bx_{i}(t)$ and  $\bg_{i}(t)$ ($\bx^{j}(t)$ and  $\bg^{j}(t)$) to show $i$th row ($j$th column) of $\bx(t)$ and $\bg(t)$, respectively. Although, $\bx_{i}(t)$ and  $\bg_{i}(t)$ are row vectors, when we use them individually, for notational convenience, we consider them column vectors.

In this paper, we are dealing with the dynamics of the $n$ agents estimates of an optimizer $z^*\in \mathcal{Z}$ which we denote them by $\bx_i(t)$ for all $i\in[n]$. Therefore, we view $\bx(t)$ as a \textbf{vector} of $n$ elements in the vector space $\R^m$. 

A distributed solution of \eqref{eqn:mainProblem} was first proposed in \cite{nedic2009distributed} using the following deterministic dynamics
\begin{align*}
    {\bx}_i(t+1)&=\sum_{j=1}^n w_{ij}(t+1){\bx}_j(t)-\alpha(t){\bg}_i(t)
\end{align*}
for $t\geq t_0$ for an initial time $t_0\in \N$, initial conditions $\bx_i(t_0)\in \R^m$ for all $i\in [n]$, where $\bg_i(t)\in\Rd$ is a subgradient of $f_{i}(z)$ at $z=\bx_{i}(t)$ for $i\in[n]$, and $\{\alpha(t)\}$ is a step-size sequence (in \cite{nedic2009distributed} the constant step-sizes variation of this dynamics was studied). {We simply refer to this dynamics as the averaging-based distributed optimization solver}. We can compactly write the above dynamics as
\begin{align}\label{eqn:MainLineAlgorithm}
    {\bx}(t+1)&=W(t+1){\bx}(t)-\alpha(t){\bg}(t),
\end{align}
where, ${\bg}(t)=[\bg_1(t),\ldots,\bg_n(t)]^T$ is the vector of the sub-gradient vectors and matrix multiplication should be understood over the vector-field $\Rd$, i.e., 
\[[W(t+1){\bx}(t)]_i\triangleq\sum_{j=1}^nw_{ij}(t+1)\bx_j(t).\] 
%In fact, $i$th row of $\bx(t)$ is the values of node $i$ in time $t$. 
In distributed optimization, the goal is to find distributed dynamics $\bx_i(t)$s such  that $\lim_{t\to\infty}\bx_i(t)=z$ where $z\in \mathcal{Z}$ for all $i\in[n]$.

\subsection{Our Contribution}
 In the paper, we consider the random variation of \eqref{eqn:MainLineAlgorithm}, i.e., when $\{W(t)\}$ is a chain of random matrices. This random variation was first studied in \cite{lobel2010distributed} where to ensure the convergence, it was assumed that this sequence is {doubly stochastic almost surely} and i.i.d.. This was generalized to random networks that is Markovian in \cite{lobel2011distributed}. The dynamics~\eqref{eqn:MainLineAlgorithm} with i.i.d.\ weight matrices that are row-stochastic almost surely and column-stochastic \textit{in-expectation} was studied in \cite{7849184}. A special case of \cite{7849184} is the asynchronous gossip algorithm that were introduced in \cite{5585721}. In this work, we provide an overarching framework for the study of~\eqref{eqn:MainLineAlgorithm} with random weight matrices that are row-stochastic almost surely and column-stochastic in-expectation, that are not even independent in general. More precisely, we have following assumptions.
\begin{assumption}[Stochastic Assumption]\label{assump:AssumpOnStochasticity}
We assume that the weight matrix sequence $\{W(t)\}$, adapted to a filtration $\{\F(t)\}$, satisfies
\begin{enumerate}[(a)]
    % \item $\{W(t)\}$ is an independent random matrix sequence.
    \item For all $t\geq t_0$, $W(t)$ is row-stochastic almost surely.
    \item For every $t> t_0$, $\E[W(t)\mid \F(t-1)]$ is column-stochastic (and hence, doubly stochastic) almost surely.
\end{enumerate}
\end{assumption}

%The problem with being doubly stochastic arises when we do not have $W(t)$ explicitly. More specifically,  we just know whether two nodes $i,j$ are connected without knowing $w_{ij}(t)$. To overcome this problem, in \cite{nedic2014distributed}, the subgradient-push method is proposed. In the subgradient-push method, $\{W(t)\}$ is constructed to be column-stochastic as follows. Every node divide its value by the number of its neighbor, and in this way, $\{W(t)\}$ is column-stochastic. Therefore, this method is sensitive to the edge failure because when a transmission fail, $W(t)$ is not column-stochastic anymore. 

%In this paper, we sanother solution is provided that need a weaker condition that double-stochasticity and is more sustainable in failure. To have consensus, we still keep being row-stochastic, but we use being column-stochastic on average. More precisely, we consider the following assumption.Almost all the \textit{simple} averaging-based distributed optimization algorithms over deterministic or random networks \cite{} rely on such an assumption. %

Similar to other works in this domain, our goal is to ensure that $\lim_{t\to\infty}\bx_i(t)=z$ {\it almost surely} for some optimal $z\in \mathcal{Z}$ for all $i\in[n]$. To reach such a consensus value, we need to ensure enough flow of information between the agents, i.e., the associated graph sequence of $\{W(t)\}$ satisfies some form of connectivity over time. More precisely, we assume the following connectivity conditions. 
\begin{assumption}[Conditional $B$-Connectivity Assumption]\label{assump:AssumpOnConnectivity}
We assume that for all $t \geq t_0$
\begin{enumerate}[(a)]
    \item \label{cond:connecta} Every node in $\G(W(t))$ has a self-loop, almost surely.
    \item \label{cond:bconnectivity}There exists an integer $B>0$ such that the random graph $\G_B(t)=([n],\EE_{B}(t))$ where
    \begin{align*}
       \EE_{B}(t)=\bigcup_{\tau=tB+1}^{(t+1)B}\EE(\E[W(\tau)|\F(tB)])
    \end{align*}
    has a spanning rooted tree almost surely.
\end{enumerate}
\end{assumption}

Note that Assumption~\ref{assump:AssumpOnConnectivity}-\eqref{cond:bconnectivity} is satisfied if the random graph with vertex set $[n]$ and the edge set 
% \begin{align*}
    $\bigcup_{\tau=tB+1}^{(t+1)B}\EE(\E[W(\tau)|\F(\tau-1)])$
% \end{align*}
has a spanning rooted tree almost surely.
% Most papers in distributed optimization assume that $W(t)$ is both row-stochastic and column-stochastic (that is double-stochastic). The problem with being double-stochastic appears when we do not have $W(t)$ explicitly. More specifically,  we just know whether two nodes $i,j$ are connected without knowing $w_{ij}(t)$. To overcome this problem, in \cite{nedic2014distributed}, the subgradient-push method is proposed. In the subgradient-push method, $\{W(t)\}$ is constructed to be column-stochastic as follows. Every node divide its value by the number of its neighbor, and in this way, $\{W(t)\}$ is column-stochastic. Therefore, this method is sensitive to the edge failure because when a transmission fail, $W(t)$ is not column-stochastic anymore. In this paper, another solution is provided that need a weaker condition that double-stochasticity and is more sustainable in failure. To have consensus, we still keep being row-stochastic, but we use being column-stochastic on average. More precisely, we consider the following assumption.

Finally, we assume the following standard condition on the step-size sequence $\{\alpha(t)\}$.
\begin{assumption}[Assumption on Step-size]\label{assump:AssumpOnStepSize}
For the step-size sequence $\{\alpha(t)\}$, we assume that $0<\alpha(t)\leq Kt^{-\beta}$ for some $K,\beta>0$ and all $t\geq t_0$, $\lim_{t\to\infty}\frac{\alpha(t)}{\alpha(t+1)}=1$, and 
\begin{align}\label{eqn:assumptionOnAlpha}
    \sum_{t=t_0}^{\infty}\alpha(t)=\infty~~\mbox{ and }~~\sum_{t=t_0}^{\infty}\alpha^2(t)<\infty.
\end{align}
\end{assumption}
The main result of this paper is the following theorem.
\begin{theorem}\label{thm:MainTheoremDisOpt}
Under the Assumptions \ref{assump:AssumpOnFunction}-\ref{assump:AssumpOnStepSize} on the model and the dynamics~\eqref{eqn:MainLineAlgorithm}, {$\lim_{t\to\infty}\bx_i(t)=z^*$ almost surely for all $i\in[n]$ and all initial conditions $\bx_i(t_0)\in \Rd$, where $z^*$ is a random vector that is supported on the  optimal set $\mathcal{Z}$. }
\end{theorem}

Before continuing with the technical details of the proof, let us first discuss some of the higher-level implications of this result: 

\textbf{1. Gossip-based sequential solvers}:   Gossip algorithms, which were originally studied in \cite{boyd2006randomized,aysal2009broadcast}, have been used in solving distributed optimization problems  \cite{5585721,7165636}. In gossip algorithms,  at each round, a node randomly wakes up and shares its value with all or some of its neighbors. However, it is possible to leverage Theorem \ref{thm:MainTheoremDisOpt} to synthesize algorithms that do not require choosing a node {independently and uniformly at random} or use other coordination methods to update information at every round. An example of such a scheme is as follows:
\begin{example}
Consider a connected {\it undirected} network\footnote{The graphs do not need to be time-invariant, and this example can be extended to processes over underlying time-varying graphs.} $\G=([n],E)$. Consider a token that is handed sequentially in the network and initially it is handed to an arbitrary agent $\ell(0)\in[n]$ in the network. If at time $t\geq 0$, agent $\ell(t)\in[n]$ is in the possession of the token, it chooses one of its neighbors $s(t+1)\in [n]$ randomly and by flipping a coin, i.e., with probability $\frac{1}{2}$ shares its information to $s(t+1)$ and passes the token and with probability $\frac{1}{2}$ keeps the token and asks for information from $s(t+1)$. It means 
\begin{align*}
    \ell(t+1)=\begin{cases}\ell(t),& \mbox{with probability $\frac{1}{2}$}\\
    s(t+1),& \mbox{with probability $\frac{1}{2}$}\end{cases}.
\end{align*}
Finally, the agent $\ell(t+1)$, who has the token at time $t+1$ and is receiving the information, does 
\begin{align*}
    \bx_{\ell(t+1)}(t+1)=\frac{1}{2}(\bx_{s(t+1)}(t)+\bx_{\ell(t)}(t))-\alpha(t)\bg_{\ell(t+1)}(t).
\end{align*}
For the other agents $i\not=\ell(t+1)$, we set
\begin{align*}
     \bx_i(t+1)=\bx_i(t)-\alpha(t)\bg_i(t).
\end{align*}
        % \begin{align}\label{eqn:tokenaround}
        %     \bx_{i}(t+1)&=
        %     \begin{cases} 
        %         \frac{1}{2}(\bx_{s(t)}(t)+\bx_{\ell(t)}(t))-\alpha(t)\bg_i(t)&\\
        %         &\hspace{-.7cm}\mbox{if $i=\ell({t}+1)$}\\
        %         \bx_i(t)-\alpha(t)\bg_i(t)&\mbox{otherwise}.
        %     \end{cases}
        % \end{align}

Let $\F(t)=\sigma(\bx(0),\ldots,\bx(t),\ell(t))$, and the weight matrix $W(t)=[w_{ij}(t)]$ be
\begin{align*}
    w_{ij}(t)=\begin{cases}\frac{1}{2},& i=j=\ell(t)\\
    \frac{1}{2},& i=\ell(t),j\in\{s(t),\ell(t-1)\}\setminus\{\ell(t)\}\\
    1,& i=j\not=\ell(t)\\
    0,& \mbox{otherwise}
    \end{cases},
\end{align*}
which is the weight matrix of this scheme. Note that  $\E[W(t)|\F(t-1)]=V(\ell(t-1))$ where $V(h)=[v_{ij}(h)]$ with
\begin{align*}
    v_{ij}(h)=\begin{cases}\frac{3}{4},& i=j=h\\
    \frac{1}{4\delta_i},& i=h,(i,j)\in E \\
    \frac{1}{4\delta_i},& j=h,(i,j)\in E\\
    1,& i=j\not=h\\
    0,& \mbox{otherwise}
    \end{cases},
\end{align*}
where $\delta_i$ is the degree of the node $i$. Note that the matrix $\E[W(t)|\F(t-1)]$ is doubly stochastic, satisfies Assumption \ref{assump:AssumpOnConnectivity}-(\ref{cond:connecta}), and only depends on $\ell(t-1)$. Now, we need to check whether $\{W(t)\}$ satisfies Assumption \ref{assump:AssumpOnConnectivity}-(\ref{cond:bconnectivity}).   
We have
\begin{align*}
    \E[W(t+n)|\F(t)]
    % &=\E[\E[V_{s(t+n)}(t+n)|\F(t+1)]|\F(t)]\cr
    &\stackrel{}{=}\E[\E[W(t+n)\mid\F(t+n-1)]\mid\F(t)]\cr
    &\stackrel{}{=}\E[V(\ell(t+n-1))\mid \F(t)]\cr
     &=\E\left[\sum_{i=1}^nV(\ell(t+n-1))1_{\{\ell(t+n-1)=i\}}\bigg{|}\F(t)\right]\cr
%    &=\sum_{i=1}^{n}\E[W(t+n)1_{\{\ell(t+n-1)=i\}}|\F(t)]\cr
    &\stackrel{}{=}\sum_{i=1}^{n}\E[V(i)1_{\{\ell(t+n-1)=i\}}\mid\F(t)]\cr
    &\stackrel{}{=}\sum_{i=1}^{n}V(i)\E[1_{\{\ell(t+n-1)=i\}}\mid\F(t)]. 
\end{align*}
If the network is  connected, starting from any vertex, after $n-1$ steps, the probability of reaching any other vertex is at least $(2\Delta)^{-(n-1)}>0$, where $\Delta\triangleq \max_{i\in [n]}\delta_i$. Therefore, we have  $\E[1_{\{\ell(t+n-1)=i\}}|\F(t)]>0$ for all $i\in[n]$ and $t$, and hence, Assumption \ref{assump:AssumpOnConnectivity}-(\ref{cond:bconnectivity}) is satisfied with $B=n$.
\end{example}

     \textbf{2. Robustness  to link-failure}: 
     Our result shows that \eqref{eqn:MainLineAlgorithm} is robust to random link-failures. Note that the results such as \cite{lobel2010distributed} will not imply the robustness of the algorithms to link failure as it assumes that the resulting weight matrices remain doubly stochastic. To show the robustness of averaging-based solvers, suppose that we have a deterministic doubly stochastic sequence $\{A(t)\}$, and suppose that each link at any time $t$ fails with some probability $p(t)>0$. More precisely, let $B(t)$ be a failure matrix where $b_{ij}(t)=0$ if a failure on link $(i,j)$ occurs at time $t$ and otherwise $b_{ij}(t)=1$ and we have
    \begin{align}\label{eqn:E1bijtuFt1begincases}
        \E[b_{ij}(t)| \F(t-1)]=1-p(t),
    \end{align}
    for $i,j\in[n]$.
    For example, if $B(t)$ is independent and identically distributed, i.e.,
    \[b_{ij}(t)=\begin{cases} 0,&\text{ with probability  $p$}\cr 
    1,& \text{ with probability  $1-p$}\end{cases},\]
    then $B(t)$ satisfies \eqref{eqn:E1bijtuFt1begincases}.
    {Define $W(t)=[w_{ij}(t)]$ as follows
    \begin{align*}
            w_{ij}(t)\triangleq\begin{cases}
            a_{ij}(t)b_{ij}(t),&i\not=j\cr
            1-\sum_{j\not=i}a_{ij}(t)b_{ij(t)},& i=j
            \end{cases}.
        \end{align*}
        % Here, $w_{ij}(k)$ models the weights with a i.i.d. link-failure. 
        Note that $W(t)$ is row-stochastic, and since  $A(t)$ is column-stochastic, $\E[W(t)|\F(t-1)]$ is column-stochastic. Thus, Theorem \ref{thm:MainTheoremDisOpt}, using $W(t)$, translates to a theorem on robustness of the distributed dynamics \eqref{eqn:MainLineAlgorithm}: as long as the connectivity conditions of Theorem~\ref{thm:MainTheoremDisOpt} holds, the dynamics will reach a minimizer of the distributed problem almost surely. For example, if the link failure probability satisfies $p(t)\leq \bar{p}$ for all $t$ and some $\bar{p}<1$, our result implies that the result of Proposition 4 in \cite{do1} (for unconstrained case) would still hold under the above link-failure model. It is worth mentioning that if $\{A(t)\}$ is time-varying, then $\E[W(t)]$ would be time-varying and hence, the previous results on distributed optimization using i.i.d.\ row-stochastic weight matrices that are column-stochastic in-expectation \cite{7849184} would not imply such a robustness result.  
    }

\section{Autonomous Averaging Dynamics}\label{sec:ConsensusWithoutPerturbation}
To prove Theorem \ref{thm:MainTheoremDisOpt}, we need to study the time-varying distributed averaging dynamics with a particular control input (gradient-like dynamics). To do this, first we study the autonomous averaging dynamics (i.e., without any input) and then, we use the established results to study the controlled dynamics.

For this, consider the time-varying distributed averaging dynamics 
\begin{align}\label{eqn:unperturbed}
    \bx(t+1)=W(t+1)\bx(t),
\end{align}
where $\{W(t)\}$ satisfying Assumption \ref{assump:AssumpOnConnectivity}. 
Defining transition matrix
\begin{align*}
    \Phi(t,\tau)\triangleq W(t)\cdots W(\tau+1),
\end{align*}
and $\Phi(\tau,\tau)=I$, we have $\bx(t)=\Phi(t,\tau)\bx(\tau)$. {Note that since $W(t)$s are row-stochastic matrices (a.s.) and the set of row-stochastic matrices is a semi-group (with respect to multiplication), the transition matrices $\Phi(t,\tau)$ are all row-stochastic matrices (a.s.). } 

{We say that a chain $\{W(t)\}$ achieves \textit{consensus} for the initial time $t_0 \in \mathbb{N}$ if for all $i$,
$\lim_{t\rightarrow\infty }\|\bx_i(t)-\tilde{x}\|=0$ almost surely,
for all choices of initial condition $\bx(t_0)\in (\Rd)^n$ in \eqref{eqn:unperturbed} and some random vector $\tilde{x}=\tilde{x}_{\bx(t_0)}$.}
It can be shown that an equivalent condition for consensus (for time $t_0$) that  $\lim_{\tinf}\Phi(t,t_0)=e \pi^T(t_0)$ for a random stochastic vector $\pi(t_0,\omega)\in \mathbb{R}^n$, almost surely where $\omega\in \Omega$ is a sample point. 

For a matrix $A=[a_{ij}]$, let
\begin{align*}
    \diamm(A)=\max_{i,j\in[n]}\frac{1}{2}\sum_{\ell=1}^n |a_{i\ell}-a_{j\ell}|,
\end{align*}
and the mixing parameter
\begin{align*}
    \Lambda(A)=\min_{i,j\in[n]}\sum_{\ell=1}^n \min\{a_{i\ell},a_{j\ell}\}.
\end{align*}
Note that for a row-stochastic matrix $A$, $\diamm(A)\in [0,1]$.
For a vector $\bx=[\bx_{i}]$ where $\bx_i\in \R^m$ for all $i$, let
\begin{align*}
    \diam(\bx)=\max_{i,j\in[n]}\|\bx_i-\bx_j\|.
\end{align*}
% where $\|\cdot\|$ is an arbitrary underlying norm in $\Rd$. 
Note that 
$\diam(\bx)\leq 2\max_{i\in[n]}\|\bx_i\|$.
% and $\diamm(A)=\max_{j}\diam(A^j)$, where $A^{(j)}$ is the $j$th column of $A$.  
Also, 
if we have consensus, then $\lim_{t\to\infty}\diam(\bx(t))=0$ and $\lim_{t\to\infty}\diamm(\Phi(t,t_0))=0$ and in fact, the reverse implications are true \cite{chatterjee1977towards}, i.e., a chain achieves consensus (for time $t_0$) if and only if $\lim_{t\to\infty}\diam(\bx(t))=0$ for all $\bx(t_0)\in \Rdn$ or $\lim_{t\to\infty}\diamm(\Phi(t,t_0))=0$.

The following results relating the above quantities are useful for our future discussions. 
\begin{lemma}[\cite{hajnal1958weak,shen2000geometric}]\label{lem:lambdadiam}
For $n\times n$ row-stochastic matrices $A,B$, we have
\begin{align*}
    \diamm(AB)\leq(1-\Lambda(A))\diamm(B).
\end{align*}

\end{lemma}
\begin{lemma}\label{lem:DiamBehrouz}
For any $n\times n$ row-stochastic matrices $A,B$, we have
\begin{enumerate}[(a)]
    \item $\diam(A\bx)\leq\diamm(A)\diam(\bx)$ for all $\bx \in(\Rd)^n$, 
    \item $\diam(\bx+\by)\leq\diam(\bx)+\diam(\by)$ for all $\bx,\by\in (\Rd)^n$,
    \item $\diamm(A)=1-\Lambda(A)$,
    \item $\diamm(AB)\leq \diamm(A)\diamm(B)$, and 
    \item $\left\|\bx_i-\sum_{j=1}^{n}\pi_j\bx_j\right\|\leq \diam(\bx)$ for all $i\in[n]$, $\bx\in (\Rd)^n$, and any stochastic vector $\pi\in [0,1]^n$ (i.e., $\sum_{i=1}^n\pi_i=1$). 
    % \item  $\diamm(A)\leq 1$.
\end{enumerate}
\end{lemma}
\pf The proof is provided in Appendix.
\qed
 
%  We may comment that the above result holds for any underlying norm $\|\cdot\|$ over $\Rd$. 

%Note that upper bounds for this quantity can be generalized to upperbounds on $\E[\diamm(\Phi(t_1,\tau_1))]$ as by Lemma \ref{lem:lambdadiam}-(e), $\diamm(\Phi(t_2,\tau_2))\leq 1$ and hence,  
%
The main goal of this section is to obtain an exponentially decreasing upper bound (in terms of $t_1-\tau_1$ and $t_2-\tau_2$) on $\E[\diamm(\Phi(t_2,\tau_2))\diamm(\Phi(t_1,\tau_1))\mid \F(\tau_1)]$.

Using this result and the connectivity assumption \ref{assump:AssumpOnConnectivity}, we can show that the transition matrices $\Phi(t,s)$ become mixing \textit{in-expectation} for large enough $t-s$. 
\begin{lemma}\label{lem:ExpectationLambda}
Under Assumption \ref{assump:AssumpOnConnectivity} (Connectivity), there exists a parameter $\theta>0$ such that for every $s\geq t_0$, we have almost surely
\begin{align*}
    \E[\Lambda(\Phi((n^2+s)B,sB))\mid\F(sB)]\geq\theta.
\end{align*}
\end{lemma}
\pf
% Without loss of generality, assume that $s=0$.
Fix $s\geq 0$.  
Let $\mathbb{T}$ be the set of all collection of edges $E$ such that the graph $([n],E)$ has a spanning rooted tree, 
and for $k\in[n^2]$,
\begin{align*}
    \EE_{B}^{}(k)\triangleq\bigcup\limits_{\tau=(s+k-1)B+1}^{(s+k)B}\EE^{}(\E[W(\tau)|\F((s+k-1)B)]).
\end{align*}
 For notational simplicity, denote $\F(sB)$ by $\F$ and $\F((s+k)B)$ by $\F_k$ for $k\in[n^2]$. Let $V=\{\omega\mid\forall k~\EE_{B}^{}(k)\in\mathbb{T}\}$.
From Assumption \ref{assump:AssumpOnConnectivity}, we have $P(V)=1$. For $\omega\in V$ and $k\geq 1$, define the random graph $([n],\mathcal{T}_k)$ on $n$ vertices by
\begin{align*}
    \mathcal{T}_k&=\begin{cases}
    \mathcal{T}_{k-1},& \mbox{if } \mathcal{T}_{k-1}\in \mathbb{T}\\
    \mathcal{T}_{k-1}\cup\{u_k\},& \mbox{if } \mathcal{T}_{k-1}\not\in \mathbb{T}
    \end{cases},
\end{align*}
with $ \mathcal{T}_0=\emptyset$, where  
\begin{align}\label{eqn:ukinEEBkcap}
    u_k\in \EE_{B}^{}(k)\cap\overline{\mathcal{T}}_{k-1},
\end{align}
and $\overline{\mathcal{T}}_{k}$ is the edge-set of the complement graph of  $([n],{\mathcal{T}}_{k})$. Note that since $\EE_{B}^{}(k)$ has a spanning rooted tree, if $\mathcal{T}_{k-1}\not\in \mathbb{T}$,  then $\EE_B(k)$ should contain an edge that does not belong to ${\mathcal{T}}_{k-1}$, which we identify it as $u_k$ in \eqref{eqn:ukinEEBkcap}. Hence, $\mathcal{T}_k$ is well-defined.
Since there  are  at  most $n(n-1)$ potential  edges  in  a  graph  on $n$ vertices, $\mathcal{T}_{n^2}$ has a spanning rooted tree for $\omega\in V$.

For $k\in[n^2]$, let
\begin{align*}
    \ED_{B}^{}(k)\triangleq\bigcup\limits_{\tau=(s+k-1)B+1}^{(s+k)B}\EE^{\nu}(W(\tau)),
\end{align*}
for some fixed $0<\nu<\gamma$, and \[\EH_{}^{}(k)\triangleq\cup_{\tau=1}^{k}\ED_{B}^{}(\tau).\]
Consider the sequences of events $\{U_k\}$ defined by $U_k\triangleq\left\{\omega\in V \mid
          \mathcal{T}_{k}\subset\EH_{}^{}(k)\right\}$,
for $k\geq 1$, and $U_0=V$.
% , and $U_0=V$, i.e., $U_k$ is the event that either the random graph $([n],\Tcal_{k-1})$ has a spanning rooted tree or the random graph $([n],\Tcal_{k})$ has an extra edge compared to $([n],\Tcal_{k-1})$. 
Note that if $\mathcal{T}_{k-1}\in \mathbb{T}$, then $\mathcal{T}_{k-1}\subset\EH_{}^{}(k-1)$ implies $\mathcal{T}_{k}\subset\EH_{}^{}(k)$, and if $\mathcal{T}_{k-1}\not\in \mathbb{T}$, then $\mathcal{T}_{k-1}\subset\EH_{}^{}(k-1)$  and $u_k\in\ED_{B}^{}(k)$ imply $\mathcal{T}_{k}\subset\EH_{}^{}(k)$. Hence, for $k\geq 1$
\begin{align}\label{eqn:DefinitionofUk}
    1_{\{U_k\}}\geq1_{\{U_{k-1}\}}1_{\{\mathcal{T}_{k-1}\not\in \mathbb{T}\}}1_{\{u_k\in\ED_{B}^{}(k)\}}+1_{\{U_{k-1}\}}1_{\{\mathcal{T}_{k-1}\in \mathbb{T}\}}.
\end{align}

% \begin{align}\label{eqn:DefinitionofUk}
%     1_{\{\mathcal{T}_{k}\in \mathbb{T}\}}&=1_{\{\mathcal{T}_{k-1}\not\in \mathbb{T}\}}1_{\{\EE_{B}^{}(k)\cap\overline{\mathcal{T}}_{k-1}\}}+1_{\{\mathcal{T}_{k-1}\in \mathbb{T}\}}.
% \end{align}
% Let's define $U_k=\{\mathcal{T}_{k}\in \mathbb{T}\}$ to be the event that the random graph $([n],\mathcal{T}_k)$ is 
% From Assumption \ref{assump:AssumpOnConnectivity}, there is a $t_kB<\tau_k\leq(t_k+1)B$ such that $u_k\in\EE(\E[W(\tau_k)])$. 
% Let us define the random process $S(k)$ for $0\leq k\leq n-1$ by 
% \begin{align*}
%     S(k)=S(k-1)+1_{\{1-w_{i_kj_k}(\tau_k)>1-\tilde{\gamma}\}},
% \end{align*}
% where $S(0)=0$. 
% We have 
% \begin{align}\label{eqn:E1UkF}
%     \E[1_{\{U_k\}}|\F]=\E[1_{\{U_k\}}1_{\{\mathcal{T}_{k-1}\not\in \mathbb{T}\}}|\F]+\E[1_{\{U_k\}}1_{\{\mathcal{T}_{k-1}\in \mathbb{T}\}}|\F]
% \end{align}
% From \eqref{eqn:DefinitionofUk}, 
On the other hand, from Tower rule, we have
% \begin{align*}
%     \E[1_{\{U_k\}}1_{\{\mathcal{T}_{k-1}\in \mathbb{T}\}}|\F]=\E[1_{\{U_{k-1}\}}1_{\{\mathcal{T}_{k-1}\in \mathbb{T}\}}|\F],
% \end{align*}
% and
\begin{align}\label{eqn:E1Uk1mathcalTk-1notinmathbbT}
    \E[1_{\{U_{k-1}\}}&1_{\{\mathcal{T}_{k-1}\not\in \mathbb{T}\}}1_{\{u_k\in\ED_{B}^{}(k)\}}\mid \F]\\
     &=\E[1_{\{U_{k-1}\}}1_{\{\mathcal{T}_{k-1}\not\in \mathbb{T}\}}\E[1_{\{u_k\in\ED_{B}^{}(k)\}}|\F_{k-1}]\mid\F].\nonumber
\end{align}
 Let $u_k(\omega)=(j_k(\omega),i_k(\omega))$.
 Since $u_k\in \EE_{B}(k)$, there exists $(s+k-1)B<\tau_k\leq (s+k)B$ such that 
 \[u_k\in\EE^{}(\E[W(\tau_k)|\F_{k-1}]),\]
 and, we have
\begin{align}\label{eqn:E1nu11wikjkkgeq}
      \E[(1-\nu) 1_{\{1-w_{i_kj_k}(\tau_k)\geq 1-\nu\}}|\F_{k-1}]\nonumber 
      &\leq \E[1-w_{i_kj_k}(\tau_k)|\F_{k-1}]
      \leq  1-\gamma.
\end{align}
Therefore, 
\begin{align*}
% \label{eqn:E[1_{{1a(tau)>1tilde{gamma}}}]}
    \E[1_{\{u_k\in\ED_{B}^{}(k)\}}\mid\F_{k-1}]&\geq\E[1_{\{u_k\in\EE^{\nu}(W(\tau_k))\}}\mid\F_{k-1}]\cr
    &=\E[ 1_{\{w_{i_kj_k}(\tau_k)> \nu\}}\mid\F_{k-1}]\cr
    &=1-\E[ 1_{\{w_{i_kj_k}(\tau_k)\leq \nu\}}\mid\F_{k-1}]\cr
    &=1-\E[ 1_{\{1-w_{i_kj_k}(\tau_k)\geq 1-\nu\}}\mid\F_{k-1}]\cr
    &\geq 1-\frac{1-\gamma}{1-\nu}\triangleq{p}>0,
\end{align*}
which holds as $\nu<\gamma$.
% We have 
% \begin{align*}
%     \E[1_{\{S(n-1)\geq 1\}}\mid \F]&\leq \E[S(n-1)\mid \F]\cr
%     &=\sum_{k=1}^{n-1}\E[1_{\{1-w_{i_kj_k}(\tau_k)>1-\tilde{\gamma}\}}\mid \F]\cr
%     &\leq (n-1)\frac{1-\gamma}{1-\tilde{\gamma}},
% \end{align*}
% which follows from \eqref{eqn:E[1_{{1a(tau)>1tilde{gamma}}}]}.
% \begin{align*}
%     \E[1_{\{U\}}\mid \F]&\geq \E[1_{\{S(n-1)=0\}}\mid \F]\cr
%     &=1-\E[1_{\{S(n-1)\geq 1\}}\mid \F]\cr 
%     &\geq 1-(n-1)\frac{1-\gamma}{1-\tilde{\gamma}}\triangleq\tilde{p}>0.
% \end{align*}
This inequality and \eqref{eqn:E1Uk1mathcalTk-1notinmathbbT} imply that
\begin{align*}
    \E[1_{\{U_{k-1}\}}1_{\{\mathcal{T}_{k-1}\not\in \mathbb{T}\}}&1_{\{u_k\in\ED_{B}^{}(k)\}}\mid \F]\geq
     p\E[1_{\{U_{k-1}\}}1_{\{\mathcal{T}_{k-1}\not\in \mathbb{T}\}}\mid\F].
\end{align*}
Therefore, \eqref{eqn:DefinitionofUk} implies 
\begin{align*}
    \E[1_{\{U_k\}}|\F]
    &\geq p\E[1_{\{U_{k-1}\}}1_{\{\mathcal{T}_{k-1}\not\in \mathbb{T}\}}|\F]+\E[1_{\{U_{k-1}\}}1_{\{\mathcal{T}_{k-1}\in \mathbb{T}\}}|\F]\cr
    &\geq p\left(\E[1_{\{U_{k-1}\}}1_{\{\mathcal{T}_{k-1}\not\in \mathbb{T}\}}|\F]+\E[1_{\{U_{k-1}\}}1_{\{\mathcal{T}_{k-1}\in \mathbb{T}\}}|\F]\right)\cr
    &=p\E[1_{\{U_{k-1}\}}|\F],
\end{align*}
and hence $\E[1_{\{U_k\}}|\F]\geq p^k$. 
Finally, since $\mathcal{T}_{n^2}$ has a spanning rooted tree, from Lemma~1 in \cite{lobel2010distributed}, we have  
\begin{align*}
    \Lambda(W(n^2B,\omega)\cdots W(n(n-1)B+n-1,\omega)\cdots W(1,\omega))\geq \nu^{n^2B},
\end{align*}
for $\omega\in U_{n^2}$. Therefore, we have 
\begin{align*}
    \E[\Lambda(\Phi((n^2+s)B,sB))\mid\F(sB)]&\geq \nu^{n^2B}\E[1_{\{U_{n^2}\}}\mid \F]\geq \nu^{n^2B}{p}^{n^2}\triangleq\theta>0, 
\end{align*}
which completes the proof.\qed
Finally, we need the following result, which is proved in Appendix, to prove the main result of this section.
\begin{lemma}\label{lem:ProductwithHistory}
For a random process $\{Y(k)\}$, adapted to a filtration $\{\F(k)\}$, let 
\[\E[Y(k)|\F(k-1)]\leq a(k)\]
for $K_1\leq k\leq  K_2$ almost surely, where $K_1\leq K_2$ are arbitrary positive integers and $a(k)$s are (deterministic) scalars. Then, we have almost surely
\begin{align*}
    \E\left[\prod_{k=K_1}^{K_2}Y(k)\bigg{|}\F(K_1-1)\right]\leq \prod_{k=K_1}^{K_2}a(k).
\end{align*}
\end{lemma}

Now, we are ready to prove the main result for the convergence rate of the autonomous random averaging dynamics. 
\begin{lemma}\label{lem:ExpectationDiam}
Under Assumption \ref{assump:AssumpOnConnectivity} (Connectivity), there exist $0<C$ and $0\leq \lambda<1$ such that for every $t_0\leq \tau_1\leq t_1$ and $t_0\leq \tau_2\leq t_2$ with $\tau_1\leq\tau_2$, we have almost surely 
\begin{align*}
    \E[\diamm(\Phi(t_2,\tau_2))\diamm(\Phi(t_1,\tau_1))|\F(\tau_1)]\leq C\lambda^{t_1-\tau_1}\lambda^{t_2-\tau_2}.
\end{align*}
\end{lemma}
\pf First,  we prove 
\begin{align}\label{eqn:EdiammPhittauFtau}
    \E[\diamm(\Phi(t,\tau))|\F(\tau)]\leq \tilde{C}\tilde{\lambda}^{t-\tau},
\end{align}
for some $0<\tilde{C}$ and $0\leq \tilde{\lambda}<1$. Let $s\triangleq\lceil \frac{\tau}{B}\rceil$ and $K\triangleq\lfloor\frac{t-sB}{n^2B}\rfloor$. We have 
\begin{align*}
    \E[\diamm(\Phi(t,\tau))|\F(\tau)]
    &= \E\Bigg{[}\diamm\Bigg{(}\Phi(t,sB+Kn^2B)\left[\prod_{k=1}^{K}\!\Phi(sB+kn^2B,sB+(k-1)n^2B)\!\right]\!\Phi(sB,\tau)\!\Bigg{)}\Bigg{|}\F(\tau)\Bigg{]}\cr    
    &\stackrel{(a)}{\leq} \!\E\left[\prod_{k=1}^{K}\!(1-\Lambda(\Phi(sB+kn^2B,sB+(k-1)n^2B))\Bigg{|}\F(\tau)\right]\cr
    % &\stackrel{(b)}{=} \prod_{k=1}^{K}\E\left[(1-\Lambda(\Phi(sB+kn^2B,sB+(k-1)n^2B))\right]\cr
    &\stackrel{(b)}{\leq} (1-\theta)^{K}\cr
    &\leq \tilde{C}(1-\theta)^{\frac{t-\tau}{n^2B}},
\end{align*}
where $\tilde{C}=(1-\theta)^{-1-\frac{1}{n^2}}$ and $(a)$ follows from Lemma~\ref{lem:lambdadiam}, Lemma~\ref{lem:DiamBehrouz}-(d), and {$\diamm(\Phi(.,.))\leq1$}, and $(b)$ follows from Lemma \ref{lem:ExpectationLambda} and \ref{lem:ProductwithHistory}. 
Since $\theta>0$, we have $\tilde{\lambda}\triangleq (1-\theta)^{\frac{1}{n^2B}}<1$.

To prove the main statement, we consider two cases: 
\begin{enumerate}[(i)]
    \item intervals $(\tau_1,t_1]$ and $(\tau_2,t_2]$ do not have an intersection, and 
    \item $(\tau_1,t_1]$ and $(\tau_2,t_2]$ intersect.
\end{enumerate}
For case (i), since the two intervals do not overlap, we have $t_1\leq\tau_2$, and hence
\begin{align*}
    \E[\diamm(\Phi(t_2,\tau_2))\diamm(\Phi(t_1,\tau_1))|\F(\tau_1)]
    &=\E\bigg{[}\E[\diamm(\Phi(t_2,\tau_2))|\F(\tau_2)]\diamm(\Phi(t_1,\tau_1))\bigg{|}\F(\tau_1)\bigg{]}\cr
    &\leq \tilde{C}\tilde{\lambda}^{t_1-\tau_1}\tilde{C}\tilde{\lambda}^{t_2-\tau_2},
\end{align*}
which follows from \eqref{eqn:EdiammPhittauFtau}.
For case (ii), let's write the union of the intervals $(\tau_1,t_1]$ and $(\tau_2,t_2]$ as disjoint union of three intervals:
\[(\tau_1,t_1]\cup(\tau_2,t_2]=(s_1,s_2]\cup (s_2,s_3]\cup (s_3,s_4],\]
for $s_1\leq s_2\leq s_3$ where $(s_2,s_3]\triangleq(\tau_1,t_1]\cap(\tau_2,t_2]$, $(s_1,s_2]\cup(s_3,s_4]\triangleq(\tau_1,t_1]\triangle(\tau_2,t_2]$.
% , and $(s_3,\ell_3]\triangleq(\tau_2,t_2]\setminus(\tau_1,t_1]$.
{Using this, it can be verified that}
\begin{align*}
    \E[\diamm(\Phi(t_2,\tau_2))\diamm(\Phi(t_1,\tau_1))|\F(\tau_1)]
    &\stackrel{(a)}{\leq}\!\!\E[\diamm(\Phi(s_4,s_3))\diamm^2(\Phi(s_3,s_2))\diamm(\Phi(s_2,s_1))|\F(\tau_1)]\cr
    &\stackrel{(b)}{\leq}\!\!\E[\diamm(\Phi(s_4,s_3))\diamm(\Phi(s_3,s_2))\diamm(\Phi(s_2,s_1))|\F(\tau_1)]\cr
    &\stackrel{(c)}{\leq} \tilde{C}\tilde{\lambda}^{s_2-s_1}\tilde{C}\tilde{\lambda}^{s_3-s_2}\tilde{C}\tilde{\lambda}^{s_4-s_3}\cr
    &= \tilde{C}\tilde{\lambda}^{s_2-s_1}2\tilde{C}\sqrt{\tilde{\lambda}}^{2(s_3-s_2)}\tilde{C}\tilde{\lambda}^{s_4-s_3}\cr
    &\leq \tilde{C}^3\sqrt{\tilde{\lambda}}^{t_1-\tau_1}\sqrt{\tilde{\lambda}}^{t_2-\tau_2},
\end{align*}
{where $(a)$ follows from Lemma \ref{lem:DiamBehrouz}-(d) ,  $(b)$ follows from {$\diamm(A)\leq1$} for all row-stochastic matrices $A$, and $(c)$ follows from \eqref{eqn:EdiammPhittauFtau} and Lemma \ref{lem:ProductwithHistory}}. Letting $C\triangleq \max\{\tilde{C}^2,
\tilde{C}^3\}$ and $\lambda\triangleq \sqrt{\tilde{\lambda}}$, we arrive at the conclusion. 
\qed
\section{Averaging Dynamics with Gradient-Flow Like Feedback }\label{sec:ConsensusWithPerturbation}
In this section, we study the controlled linear time-varying dynamics 
\begin{align}\label{eqn:EpsilonAlgorithm}
    {\bx}(t+1)&=W(t+1){\bx}(t)+\u(t).
\end{align}
Note that the feedback $\u(t)=-\alpha(t) \bg(t)$ leads to the dynamics \eqref{eqn:MainLineAlgorithm}.  
The goal of this section is to establish bounds on the convergence-rate of $\diam(\bx)$ (to zero) in-expectation and almost surely for a class of regularized input $\u(t)$.

We start with the following two lemmas.

\begin{lemma}\label{lem:d_x_t}
For dynamics \eqref{eqn:EpsilonAlgorithm} and every $t_0\leq \tau\leq t$, we have 
\begin{align}\label{lem:lemma5}
    \diam(\bx(t))&
    \stackrel{}{\leq}\diamm(\Phi(t,\tau))\diam(\bx(\tau))+\sum_{s=\tau}^{t-1}\diamm(\Phi(t,s+1))\diam(\u(s)).
\end{align}
\end{lemma}

\pf Note that the general solution for the dynamics \eqref{eqn:EpsilonAlgorithm} is given by
\begin{align}\label{eqn:PhiPerturbation}
    \bx(t)=\Phi(t,\tau)\bx(\tau)+\sum_{s=\tau}^{t-1}\Phi(t,s+1)\u(s).
\end{align}
 Therefore, using the sub-linearity property of $\diam(\cdot)$ (Lemma \ref{lem:DiamBehrouz}-(b)), we have
\begin{align*}
% \label{eqn:afterPhiPerturbation}
    \diam(\bx(t))&{\leq}\diam(\Phi(t,\tau)\bx(\tau))+\sum_{s=\tau}^{t-1}\diam(\Phi(t,s+1)\u(s))\cr
    &{\leq}\diamm(\Phi(t,\tau))\diam(\bx(\tau))+\sum_{s=\tau}^{t-1}\diamm(\Phi(t,s+1))\diam(\u(s)),
\end{align*}
where the last inequality follows from Lemma \ref{lem:DiamBehrouz}-(a). 
\qed

\begin{lemma}\label{lem:approxfunction}
Let $\{\beta(t)\}$ be a positive (scalar) sequence such that $\lim_{t\to\infty}\frac{\beta(t)}{\beta(t+1)}=1$. Then for any $\theta\in[0,1)$, there exists some $M>0$ such that
\begin{align*}
    \sum_{s=\tau}^{t-1}\beta(s)\theta^{t-s}\leq M\beta(t),
\end{align*}
 for  all $t\geq \tau\geq t_0$.
\end{lemma}
\pf
The proof is provided in Appendix.\qed

% To analyze \eqref{eqn:EpsilonAlgorithm}, we need to understand the dynamics when $\u(t)$ is a possibly feedback mechanism, which turns out to be a random policy. For this, we have the  following technical assumption on $\u(t)$.
% \begin{assumption}\label{assum:adaptedu}
% For an independent sequence $\{W(t)\}$, we say that a (possibly random) $\u(t)$ is adapted to $\{W(t)\}$ if $\u(t)$ is measurable with respect to $\F(t)=\sigma(W(1),\ldots,W(t-1))=\sigma(\bx(0),\ldots,\bx(t))$ for all $t\geq 0$. 
% \end{assumption}

To prove the main theorem, we need to study how fast $\E[\diam(\bx(t))]$ and $\E[\diam^2(\bx(t))]$ approach to zero when the diameter of the control input $\diam(\u(t))$ goes to zero. Since $\E[\diam^2(\bx(t))]\geq \E^2[\diam(\bx(t))]$, it suffice to study  convergence rate of $\E[\diam^2(\bx(t))]$.

%%%%%%%%%%%%%%%%%%%%%%%%%%%%%%%%%%%%%%%%%%%%%%%%%%%%%%%%%%%%%%%%%%%%%%%%%%%%
\begin{lemma}\label{lem:ConsensusConveregnseVariance}
Under Assumptions  \ref{assump:AssumpOnStochasticity}, \ref{assump:AssumpOnConnectivity}, and \ref{assump:AssumpOnStepSize}, if almost surely $\diam(\u(t))<q\alpha(t)$ for some $q>0$, then we have, 
 \[\frac{\E[\diam^2(\bx(t))]}{\alpha^2(t)}\leq \hat{M}\] 
 for some $\hat{M}>0$ and all $t\geq t_0$.
\end{lemma}
\pf 
%Assuming $\lim_{\tau\to\infty}\frac{\mu(\tau+1)}{\mu(\tau)}=1$ would imply that $\mu(\tau)\not=0$ for all $\tau>\tau_0$ and some $\tau_0\geq 0$. For the rest of the proof, assume $\tau_0<\tau$. 
Taking the square of both sides of \eqref{lem:lemma5},  for  $t>\tau\geq t_0$, we have 
\begin{align*}
    \diam^2(\bx(t))
    &\stackrel{}{\leq}\diamm^2(\Phi(t,\tau))\diam^2(\bx(\tau))
    % &\quad+\sum_{s=\tau}^{t-1}\diamm^2(\Phi(t,s+1))\diam^2(u(s))\cr
    +2\diamm(\Phi(t,\tau))\diam(\bx(\tau))\sum_{s=\tau}^{t-1}\diamm(\Phi(t,s+1))\diam(\u(s))\cr
    &\qquad\qquad+\sum_{s=\tau}^{t-1}\sum_{\ell=\tau}^{t-1}\diamm(\Phi(t,s+1))\diam(\u(s))\diamm(\Phi(t,\ell+1))\diam(\u(\ell)).
\end{align*}
Taking the expectation of both sides of the above inequality, and using $\diam(\u(t))<q\alpha(t)$ almost surely, we have
\begin{align*}
    \E[\diam^2(\bx(t))]
    &\stackrel{}{\leq}\E[\diamm^2(\Phi(t,\tau))\diam^2(\bx(\tau))]
    +2\sum_{s=\tau}^{t-1}\E\big{[}\diamm(\Phi(t,\tau))\diam(\bx(\tau))\diamm(\Phi(t,s+1))\diam(\u(s))\big{]}\cr
    &\qquad\qquad+\sum_{s=\tau}^{t-1}\sum_{\ell=\tau}^{t-1}\E\big{[}\diamm(\Phi(t,s+1))\diam(\u(s))\diamm(\Phi(t,\ell+1))\diam(\u(\ell))\big{]}\cr
    &\stackrel{}{\leq}\E\bigg{[}\E[\diamm^2(\Phi(t,\tau))|\F(\tau)]\diam^2(\bx(\tau))\bigg{]}+2\sum_{s=\tau}^{t-1}\E\bigg{[}\E[\diamm(\Phi(t,\tau))\diamm(\Phi(t,s+1))|\F(\tau)]\diam(\bx(\tau))\bigg{]}\alpha(s)q\cr
    &\qquad\qquad+\sum_{s=\tau}^{t-1}\sum_{\ell=\tau}^{t-1}\E[\diamm(\Phi(t,s+1))\diamm(\Phi(t,\ell+1))]\alpha(s)\alpha(\ell)q^2.
\end{align*}
Therefore, from Lemma \ref{lem:ExpectationDiam}, we have
\begin{align*}
    \E[\diam^2(\bx(t))]
    &\stackrel{}{\leq}C\lambda^{2(t-\tau)}\E[\diam^2(\bx(\tau))]+\frac{2Cq}{\lambda}\lambda^{t-\tau}\E[\diam(\bx(\tau))]\sum_{s=\tau}^{t-1}\lambda^{t-s}\alpha(s)+\frac{Cq^2}{\lambda^2}\sum_{s=\tau}^{t-1}\sum_{\ell=\tau}^{t-1}\alpha(s)\alpha(\ell)\lambda^{t-s}\lambda^{t-\ell}\cr
    &\stackrel{}{\leq}C\lambda^{2(t-\tau)}\E[\diam^2(\bx(\tau))]+\frac{2CqM}{\lambda}\lambda^{t-\tau}\E[\diam(\bx(\tau))]\alpha(t)+\frac{Cq^2M^2}{\lambda^2}\alpha^2(t),
\end{align*}
where the last inequality follows from Lemma \ref{lem:approxfunction} and the fact that
\begin{align*}
    \sum_{s=\tau}^{t-1}\sum_{\ell=\tau}^{t-1}\alpha(s)\alpha(\ell)\lambda^{t-s}\lambda^{t-\ell}=\left(\sum_{s=\tau}^{t-1}\alpha(s)\lambda^{t-s}\right)^2.
\end{align*}
Dividing both sides of the above inequality by $\alpha^2(t)$ and noting $\frac{\alpha(\tau)}{\alpha(t)}\lambda^{t-\tau}=\prod_{\kappa=\tau}^{t-1}\frac{\alpha(\kappa)}{\alpha(\kappa+1)}\lambda$,
we have
\begin{align*}
    &\frac{\E[\diam^2(\bx(t))]}{\alpha^2(t)}
    \stackrel{}{\leq}C\frac{\E[\diam^2(\bx(\tau))]}{\alpha^2(\tau)}\left(\prod_{\kappa=\tau}^{t-1}\frac{\alpha(\kappa)}{\alpha(\kappa+1)}\lambda\right)^2+2\frac{CqM}{\lambda}\frac{\E[\diam(\bx(\tau))]}{\alpha(\tau)}\left(\prod_{\kappa=\tau}^{t-1}\frac{\alpha(\kappa)}{\alpha(\kappa+1)}\lambda\right)+\frac{Cq^2M^2}{\lambda^2}.
\end{align*}
Since $\lim_{\tau\to\infty}\frac{\alpha(\tau)}{\alpha(\tau+1)}=1$, for any $\hat{\lambda}\in (\lambda,1)$, there exists $\hat{\tau}$ such that for $\tau\geq \hat{\tau}$, we have $\frac{\mu(\tau)}{\mu(\tau+1)}\lambda\leq \hat{\lambda}$. Therefore,
\begin{align*}
    \frac{\E[\diam^2(\bx(t))]}{\alpha^2(t)}
    &\stackrel{}{\leq}C\frac{\E[\diam^2(\bx(\tau))]}{\alpha^2(\tau)}\hat{\lambda}^{2(t-\tau)}+\frac{2CqM}{\lambda}\frac{\E[\diam(\bx(\tau))]}{\alpha(\tau)}\hat{\lambda}^{t-\tau}+\frac{Cq^2M^2}{\lambda^2}.
\end{align*}
Taking the limit of the above inequality, we get 
\begin{align*}
    \limsup_{t\to\infty}\frac{\E[\diam^2(\bx(t))]}{\alpha^2(t)}
    &\stackrel{}{\leq}\lim_{t\to\infty}C\frac{\E[\diam^2(\bx(\tau))]}{\alpha^2(\tau)}\hat{\lambda}^{2(t-\tau)}+\lim_{t\to\infty}\frac{2CqM}{\lambda}\frac{\E[\diam(\bx(\tau))]}{\alpha(\tau)}\hat{\lambda}^{t-\tau}+\frac{Cq^2M^2}{\lambda^2}\cr
    &=\frac{Cq^2M^2}{\lambda^2}.
\end{align*} 
As a result, there exists an $\hat{M}>0$ such that $\frac{\E[\diam^2(\bx(t))]}{\alpha^2(t)}\leq \hat{M}$.
\qed

To prove the main theorem, we also need to show that $\diam(\bx(t))$ converges to zero almost surely (as will be proved in Lemma~\ref{lem:ConsensusConveregnse}). To do so, using the previous results, we will show (in Lemma~\ref{lem:ConsensusConveregnse}) that 
\begin{align*}
    \E[\diam(\bx(t))|\F(\tau)]\leq\E[\diamm(\Phi(t,\tau))|\F(\tau)]\diam(\bx(\tau))+K\tau^{-\gamma},
\end{align*}
for some $K,\gamma>0$. However, since $\sum_{\tau=t_0}^\infty \tau^{-\gamma}$ is not necessarily summable, we cannot use the standard Siegmund-Robbins Theorem \cite{ROBBINS1971233} to argue $\diam(\bx(t))\to 0$ based on this inequality. We will use the facts that $\E[\diamm(\Phi(t,\tau))|\F(\tau)]<1$  if $t-\tau$ is large enough, and $\diamm(\Phi(t,\tau))\leq 1$ for all $t\geq \tau \geq t$. This leads us to prove a 
martingale-type result in Lemma \ref{lem:Dooblike}. To prove this lemma, we first show the following. %\todo{add a description}

\begin{lemma}\label{lem:SLLNlike}
Consider a non-negative random process $a(t)$ such that 
\[\E[a(t+1)\mid \F(t)]\leq\tilde{\lambda}\]
for some $\tilde{\lambda}<1$ and all $t\geq 0$. For $\lambda$ satisfying $\tilde{\lambda}<{\lambda}<1$, define the sequence of stopping-times $\{t_s\}_{s\geq 0}$ by  \[t_s\triangleq\inf\{t>t_{s-1}|a(t)\leq \lambda\},\]
with $t_0=0$. Then, $\lim_{s\to\infty}(t_{s+1}-t_s)t_s^{-\beta}=0$ almost surely for all $\beta>0$. 
\end{lemma}
\pf  Let us define the martingale $S(t)$ by 
\begin{align*}
    S(t)=S(t-1)+\left(1_{\{a(t)> \lambda\}}-\E[1_{\{a(t)>\lambda\}}|\F(t-1)]\right),
\end{align*}
where $S(0)=0$. 
 Noting $|S(t+1)-S(t)|\leq 1$,  from Azuma's inequality, we have
\begin{align}\label{eqn:AzumaForSSigam}
    P(S(t+\sigma)-S(t)>\sigma\rho)\leq \exp\left(-\frac{\sigma^2\rho^2}{2\sigma}\right),
\end{align}
for all $\sigma\in \mathbb{N}$ and $\rho\in(0,1)$. For $\theta>0$, let the sequences of events 
\begin{align*}
    A_{\theta}(t)\triangleq\left\{\omega\big{|}S(t+\lfloor \theta t^\beta \rfloor)-S(t)>\lfloor \theta t^\beta \rfloor\rho\right\}.
\end{align*}
From \eqref{eqn:AzumaForSSigam}, we have
\begin{align*}
    P(A_{\theta}(t))\leq\exp\left(-\frac{1}{2}(\theta t^\beta-1)\rho^2\right), 
\end{align*}
% which follows from the fact that for every polynomial function of $h(t)$ such that $h(t)>0$ for all $t>\tau$ for some $\tau\in\R$
% \begin{align*}
%     \lim_{t\to\infty}t^2\exp(-h(t))\leq\lim_{t\to\infty}\left(\sum_{k=0}^K \frac{h^k(t)}{t^2k!}\right)^{-1}=0,
% \end{align*}
% where $\lim_{\to\infty}\frac{h^K(t)}{t^2}=\infty$.
implying $\sum_{t=1}^{\infty}P(A_{\theta}(t))<\infty$ as $\exp\left({-t^\beta}\right)\leq \frac{M}{t^2}$ for sufficiently large $M$ (depending on $\beta$). Therefore, the Borel–Cantelli Theorem implies that $P(\{A_{\theta}(t)\quad \text{i.o.}\})=0$ for all $\theta>0$.

For $\theta>0$, let the sequences of events 
\begin{align*}
    B_{\theta}(t)\triangleq\left\{\omega\bigg{|}\frac{\hat{t}-t}{t^{\beta}}>\theta \mbox{ where } \hat{t}=\inf\{\tau>t|a(\tau)\leq \lambda\}\right\}.
\end{align*}
We show that $B_{\theta}(t)\subset A_{\theta}(t)$ for all $t,\theta$.
Fix a constant $\rho\in (0,1)$ such that  $1-\frac{{\tilde{\lambda}}}{\lambda}>\rho$. Since $\E[a(\tau)\mid \F(\tau-1)]\leq \tilde{\lambda}$, we have
\begin{align*}
      \E[\lambda 1_{\{a(\tau)\geq \lambda\}}\mid \F(\tau-1)]&\leq \E[a(\tau)\mid \F(\tau-1)]\cr
     &\leq \tilde{\lambda}<\lambda(1-\rho)
\end{align*}
and hence,
\begin{align}\label{eqn:rhoE1ataugeqlambdamidFtau}
    \E[ 1_{\{a(\tau)\geq \lambda\}}\mid \F(\tau-1)]<1-\rho.
\end{align}
Let $\sigma(t)\triangleq\lfloor \theta t^\beta \rfloor$. If $\hat{t}-t>\theta t^\beta$, then 
\begin{align*}
    S(t+\sigma(t))-S(t)&= \sigma(t)-\sum_{\tau=t+1}^{t+\sigma(t)} \E[1_{\{a(\tau)>\lambda\}}|\F(\tau-1)]\cr
    &> \sigma(t)-\sigma(t)(1-\rho)=\sigma(t)\rho,
\end{align*}
which follows from \eqref{eqn:rhoE1ataugeqlambdamidFtau}.
Therefore, we have $B_{\theta}(t)\subset A_{\theta}(t)$, and hence, $P(\{B_{\theta}(t)\quad \text{i.o.}\})=0$ for all $\theta>0$.

Finally, by contradiction, we  show that $\lim_{s\to\infty}(t_{s+1}-t_s)t_s^{-\beta}=0$. Since, if $\lim_{s\to\infty}(t_{s+1}-t_s)t_s^{-\beta}\not=0$ almost surely, then $\limsup_{s\to\infty}(t_{s+1}-t_s)t_s^{-\beta}>0$ almost surely, and hence, $P\left(\limsup_{s\to\infty}(t_{s+1}-t_s)t_s^{-\beta}>\epsilon\right)>0$ for some $\epsilon>0$. Therefore, $P(\{B_{\epsilon}(t)\quad \text{i.o.}\})>0$, which is a contradiction.
\qed

% Consider the sequence of events
% \begin{align*}
%     B_t&\triangleq\{\omega|a(t,\omega)\leq \bar{\lambda}=\frac{1+\lambda}{2}\}.
% \end{align*}
% Since,  $\E[a(t+1)\mid \F_t]\leq\lambda$ and $a(t+1)\leq 1$, we have 
% \[P(B_{t+1}\mid \F_k)\geq \frac{\lambda}{\lambda+1}\]
% Since $\{a(t)\}$ is independent, we have $P(B_k)=0$. Hence $P(\cup_{k=1}^{\infty}B_k)=0$, and hence $P(\Omega\setminus\cup_{k=1}^{\infty}B_k)=1$. 
% Thus, we only need to prove $\lim_{t\to\infty}D(t,\omega)=0$ for $\omega\in B\triangleq\Omega\setminus\cup_{k=1}^{\infty}B_k$. Fix a sample path $\omega$. 
\begin{lemma}\label{lem:Dooblike}
Suppose that $\{D(t)\}$ is a non-negative random (scalar) process such that 
 \begin{align}\label{eqn:dooblike}
 D(t+1)\leq a(t+1)D(t)+b(t),\quad \mbox{almost surely}
 \end{align}
 where $\{b(t)\}$ is a deterministic sequence and $\{a(t)\}$ is an adapted process (to $\{\F(t)\}$), such that $a(t)\in[0,1]$ and 
\[\E[a(t+1)\mid \F(t)]\leq\tilde{\lambda},\]
almost surely for some $\tilde{\lambda}<1$ and all $t\geq 0$. Then, if 
\[0\leq b(t)\leq{K}{t^{-\tilde{\beta}}}\]
for some $K,\tilde{\beta}>0$, we have $\lim_{t\to\infty}D(t)t^\beta=0$, almost surely, for all $\beta<\tilde{\beta}$. 
\end{lemma}
\pf Let $t_s\triangleq\inf\{t>t_{s-1}|a(t)\leq \lambda\}$ and $t_0=0$ for some $\tilde{\lambda}<\lambda<1$, and $c(s)\triangleq \sum_{\tau=t_s+1}^{t_{s+1}-1}b(\tau)$ .
Also, define
\begin{align*}
    A\triangleq\left\{\omega\bigg{|}\lim_{s\to\infty}\frac{t_{s+1}-t_s}{t_s^{\min\{\tilde{\beta}-\beta,1\}}}=0\right\}.
\end{align*}
Note that Lemma \ref{lem:SLLNlike} implies $P(A)=1$. 
On the other hand, using \eqref{eqn:dooblike}, we have
\begin{align*}
    D(t_{s+1})&\leq D(t_{s})\prod_{\ell=t_{s}+1}^{t_{s+1}}a(\ell)+\sum_{\tau=t_s}^{t_{s+1}-1}b(\tau)\prod_{\ell=\tau+2}^{t_{s+1}}a(\ell)\cr
    &\leq D(t_{s})\lambda+c(s),
\end{align*}
where the last inequality follows from $a(t)\in [0,1]$ and $a(t_{s+1})\leq \lambda$.
Letting $R(t)=D(t)t^{\beta}$, we have
\begin{align*}
    R(t_{s+1})\leq \left(\frac{t_{s+1}}{t_s}\right)^\beta R(t_s)\lambda+c(s)t_{s+1}^{\beta}.
\end{align*}
Note that, for $\omega\in A$, we have 
% \begin{align*}
%     \lim_{s\to\infty}\left(\frac{t_{s+1}}{t_{s}}\right)^\gamma=\lim_{s\to\infty}\left(1+\frac{t_{s+1}-t_s}{t_{s}}\right)^\gamma\leq\lim_{s\to\infty}\exp\left(\frac{\gamma(t_{s+1}-t_s)}{t_{s}}\right)=1
% \end{align*}
\begin{align}\label{eqn:lim_{stoinfty}frac{t_{s+1}}{t_{s}}}
    \lim_{s\to\infty}\frac{t_{s+1}}{t_{s}}=\lim_{s\to\infty}1+\frac{t_{s+1}-t_s}{t_{s}}=1.
\end{align}
As a result, for any $\hat{\lambda}\in (\lambda,1)$, there exists $\hat{s}$ such that for $s\geq \hat{s}$, we have $\left(\frac{t_{s+1}}{t_s}\right)^\beta\lambda\leq \hat{\lambda}$, and hence
\begin{align*}
    R(t_{s+1})\leq  R(t_s)\hat{\lambda}+c(s)t_{s+1}^{\beta}.
\end{align*}
Therefore,
\begin{align*}
    R(t_s)\leq \hat{\lambda}^{s-\hat{s}} R(t_{\hat{s}})+\sum_{\tau=\hat{s}}^{s-1}c(\tau)t_{\tau+1}^{\beta}\hat{\lambda}^{s-\tau-1}.
\end{align*}
Taking the limits of the both sides,  we have
\begin{align*}
    \limsup_{s\to\infty}R(t_s)&\leq \limsup_{s\to\infty}\hat{\lambda}^{s-\hat{s}} R(t_{\hat{s}})+\sum_{\tau=\hat{s}}^{s-1}c(\tau)t_{\tau+1}^{\beta}\hat{\lambda}^{s-\tau-1}\cr
    &=\lim_{s\to\infty}c(s)t_{s+1}^{\beta},
\end{align*}
which is implied by Lemma 3.1-(a) in \cite{ram2010distributed}.
For $\omega\in A$, we have
\begin{align}\label{eqn:lim_{stoinfty}c(s)t_{s+1}^{gamma}}
    \lim_{s\to\infty}c(s)t_{s+1}^{\beta}
    &\stackrel{(a)}{\leq} \lim_{s\to\infty}\frac{K(t_{s+1}-t_s)}{t_s^{\tilde{\beta}}}t_{s+1}^{\beta}\cr
    &= \lim_{s\to\infty}\frac{K(t_{s+1}-t_s)}{t_s^{\tilde{\beta}-\beta}}\frac{t_{s+1}^{\beta}}{t_{s}^{\beta}}\cr
    &\stackrel{(b)}{=} \lim_{s\to\infty}\frac{K(t_{s+1}-t_s)}{t_s^{\tilde{\beta}-\beta}}=0.
\end{align} 
where $(a)$ follows from $b(t)\leq Kt^{-\tilde{\beta}}$, and $(b)$ follows from \eqref{eqn:lim_{stoinfty}frac{t_{s+1}}{t_{s}}}. Therefore, $\lim_{s\to\infty}R(t_s)=0$. 
% For $\omega\in A$, we have 
% \begin{align*}
%     \lim_{s\to\infty}(t_{s+1}-t_s)\frac{t_{s+1}^{\gamma}}{t_s^{\tilde{\gamma}}}
%     =\lim_{s\to\infty}\frac{t_{s+1}-t_s}{t_{s}^{\tilde{\gamma}-\gamma}}=0.
% \end{align*}
Now for any $t>0$ with $t_s\leq t<t_{s+1}$, let $\sigma(t)=s$. By the definition of $R(t)$, we have $R(t)\leq \left(\frac{t}{\sigma(t)}\right)^\beta R(t_{\sigma(t)})+c(\sigma(t))t_{\sigma(t)+1}^{\beta}$. Therefore, $\lim_{s\to\infty}R(t_s)=0$ and Inequality \eqref{eqn:lim_{stoinfty}c(s)t_{s+1}^{gamma}} imply 
\begin{align*}
    \limsup_{t\to\infty}R(t)&\leq\lim_{t\to\infty}\left(\frac{t}{\sigma(t)}\right)^\beta R(t_{\sigma(t)})+c(\sigma(t))t_{\sigma(t)+1}^{\beta}=0,
\end{align*}
which is the desired conclusion as $R(t)=D(t)t^\beta$. 
\qed

% \begin{lemma}\label{lem:Dooblike}
% Consider an adapted random process $\{a(t)\}$ (to $\{\F(t)\}$), such that $a(t)\in[0,1]$ almost surely and 
% \[\E[a(t+1)\mid \F(t)]\leq\tilde{\lambda}\]
% for some $\tilde{\lambda}<1$ and all $t\geq 0$. If $D(t+1)\leq a(t+1)D(t)+b(t)$ where $0\leq b(t)\leq{K}{t^{-\tilde{\gamma}}}$ for some $K,\tilde{\gamma}>0$, and $D(t)\geq 0$ for all $t$, then for  $\tilde{\gamma}>{\gamma}$ we have $\lim_{t\to\infty}D(t)t^\gamma=0$, almost surely. 
% \end{lemma}
Finally, we are ready to show the almost sure convergence $\lim_{t\to\infty}\diam(\bx(t))=0$ (and more) under our connectivity assumption and a regularity condition on the input $\u(t)$ for the controlled averaging dynamics \eqref{eqn:EpsilonAlgorithm}. 
\begin{lemma}\label{lem:ConsensusConveregnse}
Suppose that  $\{W(t)\}$ satisfies Assumption  \ref{assump:AssumpOnConnectivity}. Then, if $\diam(\u(t))<qt^{-\tilde{\beta}}$ almost surely for some $q\geq 0$, we have
    $\lim_{t\to\infty} {\diam(\bx(t))}{t^{\beta}}= 0$,  {almost surely}, for $\beta<\tilde{\beta}$.
\end{lemma}
\pf  
From inequality \eqref{lem:lemma5}, we have
\begin{align}\label{eqn:diambxk}
    \diam(\bx(k))
    &\leq\diamm(\Phi(k,\tau))\diam(\bx(\tau))+\sum_{s=\tau}^{k-1}\diamm(\Phi(\tau,s+1))\diam(\u(s))\cr
    &\stackrel{(a)}{\leq}\diamm(k,\tau))\diam(\bx(\tau))
    +\sum_{s=\tau}^{k-1}qs^{-\tilde{\beta}}\cr
    &\leq\diamm(\Phi(k,\tau))\diam(\bx(\tau))
    +(k-\tau)q\tau^{-\tilde{\beta}},
\end{align}
 where $(a)$ follows from $\diamm(\Phi(.,.))\leq 1$.
Let $C>0$ and $\lambda\in [0,1)$ be the constants satisfying the statement of Lemma~\ref{lem:ExpectationDiam}. Since ${\lambda}<1$, for $T=\lceil|\frac{\log C}{\log {\lambda}}|\rceil+1$, we have $\tilde{\lambda}\triangleq C{\lambda}^T< 1$.
Then, Lemma \ref{lem:ExpectationDiam} implies that
\begin{align}\label{eqn:EdiammPhihnut1hnutFhnut}
    \E[\diamm(\Phi(T(t+1),Tt))|\F(Tt)]\leq C{\lambda}^T= \tilde{\lambda}<1.
\end{align}
Let $D(t)\triangleq \diam(\bx(Tt))$.
 From inequality \eqref{eqn:diambxk}, for $\tau=Tt$ and $k=T(t+1)$, we have
\begin{align*}
    D(t+1)
    \leq\diamm(\Phi(T(t+1),Tt))D(t)
    +T^{1-\tilde{\beta}}qt^{-\tilde{\beta}}.
\end{align*}

Taking conditional expectation of both sides  of the above inequality given $\F(Tt)$, we have
\begin{align*}
    &\E[D(t+1)|\F(Tt)]\leq\E[\diamm(\Phi(T(t+1),Tt))|\F(Tt)]D(t)
    +T^{1-\tilde{\beta}}qt^{-\tilde{\beta}}.
\end{align*}
 By letting $a(t+1)\triangleq \diamm(\Phi(T(t+1),Tt))$ and $b(t)\triangleq Tqt^{-\tilde{\beta}}$, we are in the setting of Lemma \ref{lem:Dooblike}. Therefore, by Inequality \eqref{eqn:EdiammPhihnut1hnutFhnut} and $\diamm(\Phi(T(t+1),Tt))\leq 1$, the conditions of  Lemma~\ref{lem:Dooblike}  hold, and hence,  \begin{align}\label{eqn:Dttlimit}
     \lim_{t\to\infty}D(t)t^\beta=0
 \end{align} 
 almost surely. 
 
 On the other hand, letting $\tau=T\left\lfloor \frac{k}{T}\right\rfloor$ in  \eqref{eqn:diambxk} we have
  \begin{align*}
    \diam(\bx(k))\leq D\left(\left\lfloor \frac{k}{T}\right\rfloor\right)
    +T^{1-\tilde{\beta}}q\left\lfloor \frac{k}{T}\right\rfloor^{-\tilde{\beta}}. 
 \end{align*}
 Therefore, 
 \begin{align*}
    \lim_{k\to\infty}\diam(\bx(k))k^{\beta}
    &\leq\lim_{k\to\infty} D\left(\left\lfloor \frac{k}{T}\right\rfloor\right)k^{\beta}
    +T^{1-\tilde{\beta}}q\left\lfloor \frac{k}{T}\right\rfloor^{-\tilde{\beta}}k^{\beta}\cr
    &\leq\lim_{k\to\infty} D\left(\left\lfloor \frac{k}{T}\right\rfloor\right)\left(T\left(\left\lfloor \frac{k}{T}\right\rfloor+1\right)\right)^{\beta}
    +Tq(k-T)^{-\tilde{\beta}}k^{\beta}\cr&=0,
\end{align*}
where the last equality follows from \eqref{eqn:Dttlimit} and $\tilde{\beta}>\beta$.
\qed

\section{Convergence Analysis of the Main Dynamics}\label{sec:Convergence}
Finally, in this section, we will study the main dynamics \eqref{eqn:MainLineAlgorithm}, i.e., the dynamics \eqref{eqn:EpsilonAlgorithm} with the feedback policy $\u_i(t)=-\alpha(t)\bg_i(t)$ where $\bg_i(t)\in \nabla f_i(\bx_i(t))$. Throughout this section,  we let $\bar{\bx}\triangleq \frac{1}{n}e^T\bx$  for a vector $\bx\in (\R^m)^n$,

 First, we prove an inequality (Lemma~\ref{lem:FirstLyapunovDisOpt}) which plays a key role in the proof of Theorem \ref{thm:MainTheoremDisOpt} and to do so, we make use of the following result which is proven as a part of the proof of Lemma 8 (Equation (27)) in \cite{nedic2014distributed}. \begin{lemma}[\cite{nedic2014distributed}]\label{lem:AngeliaLemma}
Under Assumption \ref{assump:AssumpOnFunction}, for all $v\in\Rd$, we have
\begin{align*}
    n\lrangle{\bar{\bg}(t),\bar{\bx}(t)-\bv}\geq F(\bar{\bx}(t))-F(\bv)-2\sum_{i=1}^nL_i\|\bx_i(t)-\bar{\bx}(t)\|.
\end{align*}
\end{lemma}
\begin{lemma}\label{lem:FirstLyapunovDisOpt}
For the dynamics \eqref{eqn:MainLineAlgorithm}, under Assumption \ref{assump:AssumpOnFunction}, for all $\bv\in \Rd$, we have 
\begin{align*}
    \E[\|\bar{\bx}(t+1)-\bv\|^2|\F(t)]
    &\leq\|\bar{\bx}(t)-\bv\|^2+\alpha^2(t)\frac{L^2}{n^2}+\sum_{i=1}^{n}\|\bx_i(t)-\bar{\bx}(t)\|^2\\
    &\qquad\qquad-\frac{2\alpha(t)}{n}(F(\bar{\bx}(t))-F(\bv))+\frac{4\alpha(t)}{n}\sum_{i=1}^nL_i\|\bx_i(t)-\bar{\bx}(t)\|.
\end{align*}
\end{lemma}
\pf 
Multiplying $\frac{1}{n}e^T$ from left to both sides of \eqref{eqn:MainLineAlgorithm}, we have 
\begin{align*}
    \bar{\bx}(t+1)&=\overline{W}(t+1){\bx}(t)-\alpha(t)\bar{\bg}(t)\\
    &=\bar{\bx}(t)-\alpha(t)\bar{\bg}(t)+\overline{W}(t+1){\bx}(t)-\bar{\bx}(t),
\end{align*}
where $\overline{W}(t)\triangleq\frac{1}{n}e^TW(t)$.
Therefore, we can write
\begin{align*}
    \|\bar{\bx}(t+1)-v\|^2
    &=\|\bar{\bx}(t)-{\bv}-\alpha(t)\bar{\bg}(t)+\overline{W}(t+1){\bx}(t)-\bar{\bx}(t)\|^2\\
    &=\|\bar{\bx}(t)-{\bv}\|^2+\|\alpha(t)\bar{\bg}(t)\|^2+\|\overline{W}(t+1){\bx}(t)-\bar{\bx}(t)\|^2-2\alpha(t)\lrangle{\bar{\bg}(t),\overline{W}(t+1){\bx}(t)-{\bv}}\\
    &\qquad\qquad+2\lrangle{\bar{\bx}(t)-{\bv},\overline{W}(t+1){\bx}(t)-\bar{\bx}(t)}.
\end{align*}
Taking conditional expectation of both sides of { the above equality given 
$\F(t)$,}
we have
\begin{align*}
    \E[\|\bar{\bx}(t+1)-{\bv}\|^2|\F(t)]
    &=\|\bar{\bx}(t)-{\bv}\|^2+\|\alpha(t)\bar{\bg}(t)\|^2+\E\left[\|\overline{W}(t+1){\bx}(t)-\bar{\bx}(t)\|^2|\F(t)\right]\\
    &\quad-2\alpha(t)\lrangle{\bar{\bg}(t),\E\left[\overline{W}(t+1){\bx}(t)-{\bv}|\F(t)\right]}+2\lrangle{(\bar{\bx}(t)-{\bv}),\E\left[\overline{W}(t+1){\bx}(t)-\bar{\bx}(t)|\F(t)\right]}\\
    &{=}\|\bar{\bx}(t)-{\bv}\|^2+\|\alpha(t)\bar{\bg}(t)\|^2+\E\left[\|\overline{W}(t+1){\bx}(t)-\bar{\bx}(t)\|^2|\F(t)\right]-2\lrangle{\alpha(t)\bar{\bg}(t),\bar{\bx}(t)-{\bv}}.
\end{align*}
The last equality follows from the assumption that, $W(t+1)$ is doubly stochastic in-expectation and hence, 
\[\E[\overline{W}(t+1)|\F(t)]=\frac{1}{n}e^T,\] 
which implies  \[\lrangle{\bar{\bg}(t),\E\left[\overline{W}(t+1){\bx}(t)-{\bv}|\F(t)\right]}=\lrangle{\alpha(t)\bar{\bg}(t),\bar{\bx}(t)-{\bv}},\]
and 
\[\E\left[\overline{W}(t+1){\bx}(t)-\bar{\bx}(t)|\F(t)\right]=0.\]
Note that $\overline{W}(t+1)$ is a stochastic vector (almost surely), therefore, due to the convexity of norm-square $\|\cdot\|^2$, we get 
\begin{align*}
\|\overline{W}(t+1){\bx}(t)-\bar{\bx}(t)\|^2&\leq \sum_{i=1}^n\overline{W}_i(t+1)\|{\bx}_i(t)-\bar{\bx}(t)\|^2\cr 
&\leq \sum_{i=1}^n\|{\bx}_i(t)-\bar{\bx}(t)\|^2,
\end{align*}
{as $\overline{W}_i(t+1)\leq 1$ for all $i\in[n]$. Therefore, }
\begin{align*}
    \E[\|\bar{\bx}(t+1)-{\bv}\|^2|\F(t)]&{=}
    \|\bar{\bx}(t)-{\bv}\|^2+\|\alpha(t)\bar{\bg}(t)\|^2+\E\left[\|\overline{W}(t+1){\bx}(t)-\bar{\bx}(t)\|^2|\F(t)\right]-2\lrangle{\alpha(t)\bar{\bg}(t),\bar{\bx}(t)-{\bv}}\cr 
    &\leq \|\bar{\bx}(t)-{\bv}\|^2+\|\alpha(t)\bar{\bg}(t)\|^2+\sum_{i=1}^{n}\|\bx_i(t)-\bar{\bx}(t)\|^2-2\lrangle{\alpha(t)\bar{\bg}(t),\bar{\bx}(t)-{\bv}}. 
\end{align*}
Finally, Lemma \ref{lem:AngeliaLemma} and the fact that 
\begin{align*}
   \|\bar{\bg}(t)\|^2=\frac{1}{n^2}\left\|\sum_{i=1}^n\bg_i(t)\right\|^2\leq \frac{1}{n^2}\left(\sum_{i=1}^n\left\|\bg_i(t)\right\|\right)^2\leq \frac{L^2}{n^2},
\end{align*}
complete the proof.\qed

{\it Proof of Theorem \ref{thm:MainTheoremDisOpt}:} 
Recall that Siegmund and Robbins Theorem \cite{ROBBINS1971233} states that for a non-negative random process $\{V(t)\}$ (adapted to the filtration $\{\F(t)\}$) satisfying  
\begin{align}\label{eqn:RS}
    \E&\left[V(t+1)|\F(t)\right]=(1+a(t))V(t)-b(t)+c(t),
\end{align}
where $a(t),b(t),c(t)\geq 0$ for all $t$, and $\sum_{t=0}^{\infty}a(t)<\infty$ and $\sum_{t=1}^{\infty}c(t)<\infty$ almost surely,  $\lim_{t\to\infty}V(x(t))$ exists and $\sum_{t=1}^{\infty}b(t)<\infty$ almost surely. In order to utilize this result and Lemma \ref{lem:FirstLyapunovDisOpt}, for all $t\geq 0$, let 
\begin{align*}
V(t)&\triangleq\|\bar{\bx}(t)-z\|^2, \cr  
a(t)&\triangleq 0, \cr
b(t)&\triangleq-\frac{2\alpha(t)}{n}(F(\bar{\bx}(t))-F(z)),\text{ and }\cr
c(t)&\triangleq\alpha^2(t)\frac{L^2}{n^2}+\sum_{i=1}^{n}\|\bx_i(t)-\bar{\bx}(t)\|^2+\frac{4\alpha(t)}{n}\sum_{i=1}^nL_i\|\bx_i(t)-\bar{\bx}(t)\|,
\end{align*}
where $z\in\mathcal{Z}$. First, note that $a(t),b(t),c(t)\geq 0$ for all $t$. To invoke the Siegmund and Robbins result \eqref{eqn:RS}, we need to to prove that  $\sum_{t=0}^{\infty}c(t)<\infty$, almost surely. Since $\sum_{t=0}^{\infty}\alpha^2(t)<\infty$, it is enough to show that  
\[\sum_{t=0}^{\infty}\|\bx_i(t)-\bar{\bx}(t)\|^2<\infty\] and 
\[\sum_{t=0}^\infty\alpha(t)\|\bx_i(t)-\bar{\bx}(t)\|<\infty,\]
almost surely, for all $i\in[n]$. From Lemma~\ref{lem:DiamBehrouz}-(e) and Lemma~\ref{lem:ConsensusConveregnseVariance}, we have 
{\[\E\left[ \frac{\|\bx_i(t)-\bar{\bx}(t)\|}{\alpha(t)}\right]\leq \frac{\E\left[\diam(\bx(t))\right]}{\alpha(t)}\leq \sqrt{\hat{M}}<\infty\]}
for some $\hat{M}>0$. Therefore, we have
    \begin{align*}
       \lim_{T\to\infty}\E\left[\sum_{t=0}^{T}\alpha(t)\|\bx_i(t)-\bar{\bx}(t)\|\right]
        &=\lim_{T\to\infty}\E\left[ \sum_{t=0}^{T}\alpha^2(t)\frac{\|\bx_i(t)-\bar{\bx}(t)\|}{\alpha(t)}\right]\cr
        &=\lim_{T\to\infty}\sum_{t=0}^{T}\alpha^2(t)\E\left[ \frac{\|\bx_i(t)-\bar{\bx}(t)\|}{\alpha(t)}\right]\cr
        &\stackrel{}{\leq} \sqrt{\hat{M}}\sum_{t=0}^{\infty}{\alpha^2(t)}<\infty. 
    \end{align*}
    {Similarly, using Lemma~\ref{lem:ConsensusConveregnseVariance}, there exists some $\hat{M}>0$ such that 
    \begin{align*}
        \frac{\E\left[\|\bx_i(t)-\bar{\bx}(t)\|^2\right]}{\alpha^2(t)}\leq \frac{\E\left[\diam^2(\bx(t))\right]}{\alpha^2(t)}\leq\hat{M},
    \end{align*}
    for all $t\geq 0$, where the first inequality follows from Lemma~\ref{lem:DiamBehrouz}~-~(e).} Therefore, 
    \begin{align*}
        \lim_{T\to\infty}\E\left[\sum_{t=0}^{T}\|\bx_i(t)-\bar{\bx}(t)\|^2\right]
        &= \lim_{T\to\infty}\E\left[\sum_{t=0}^{T}\alpha^2(t)\frac{\|\bx_i(t)-\bar{\bx}(t)\|^2}{\alpha^2(t)}\right]\cr
        &\stackrel{}{\leq} \hat{M}\sum_{t=0}^{\infty}\alpha^2(t)<\infty. 
    \end{align*}
    Therefore, using Monotone Convergence Theorem, we have 
    \begin{align*}
       \E\left[\sum_{t=0}^{\infty}\alpha(t)\|\bx_i(t)-\bar{\bx}(t)\|\right]&<\infty,\text{ and }\cr 
       \E\left[\sum_{t=0}^{\infty}\|\bx_i(t)-\bar{\bx}(t)\|^2\right]&<\infty,
    \end{align*}
    which implies $\sum_{t=0}^{\infty}\alpha(t)\|\bx_i(t)-\bar{\bx}(t)\|<\infty$ and $\sum_{t=0}^{\infty}\|\bx_i(t)-\bar{\bx}(t)\|^2<\infty$,  almost surely.

    Now that we showed that $c(t)$ is almost surely a summable sequence, Siegmund and Robbins Theorem implies that almost surely \[\lim_{t\to\infty}V(t)=\lim_{t\to\infty}\|\bar{\bx}(t)-z\|^2\text{ exists},\]
    and 
    \[\sum_{t=1}^{\infty}{\alpha(t)}(F(\bar{\bx}(t))-F(z))<\infty.\]
    % The last inequality implies that 
    % \[\liminf_{t\to\infty }F(\bar{\bx}(t))-F(z^*)=0.\]
    
    % On the other hand, since $\lim_{t\to\infty}\|\bar{\bx}(t)-z^*\|^2$, the sequence $\{\bar{\bx}(t)\}$ is a bounded sequence and for almost all sample points $\omega\in \Omega$, there exists a convergent subsequece $\{\bar{x}(t_r)\}$ to some $\bar{z}$ such that 
    % \[\lim_{r\to\infty} F(\bar{x}(t_r))=F(z^*).\]
    
    % On the other hand, since $F(\cdot)$ is a continuous function 
    % \[\lim_{r\to\infty} F(\bar{x}(t_r))=F(\bar{z}),\] 
    % and hence, $\bar{z}\in \mathcal{Z}$. 
    %Therefore, if $\liminf_{t\to\infty}$

   For $z\in \mathcal{Z}$, let's define 
   \begin{align*}
       \Omega_{z}\triangleq\left\{\omega\Bigg{|}\begin{array}{cc}
             \lim\limits_{t\to\infty}\|\bar{\bx}(t,\omega)-z\| \text{ exists,}  \\ 
            \sum\limits_{t=1}^{\infty}{\alpha(t)}(F(\bar{\bx}(t,\omega))-F^*)<\infty 
       \end{array}\right\},
   \end{align*}
   where $F^*\triangleq\min_{z\in\Rd}F(z)$. Per Sigmund-Robbins result, we know that $P(\Omega_{z})=1$. Now, let $\zd\subset \mathcal{Z}$ be a countable dense subset of $\mathcal{Z}$ and let 
   \begin{align*}
       \Omega_d\triangleq\bigcap_{z\in \zd}\Omega_z. 
   \end{align*}
    Since $\zd$ is a countable set, we have $P(\Omega_d)=1$ and for $\omega\in \Omega_d$, since $\sum_{t=1}^{\infty}{\alpha(t)}(F(\bar{\bx}(t,\omega))-F^*)<\infty $ and $\alpha(t)$ is not summable, we have 
    \[\liminf_{t\to\infty}F(\bar{\bx}(t))=F^*.\]
   This fact and the fact that $F(\cdot)$ is a continuous function
    implies that for all $\omega\in \Omega_d$, we have $\liminf_{t\to\infty}\|\bar{\bx}(t,\omega)-z^*(\omega)\|=0$ for some $z^*(\omega)\in\mathcal{Z}$.
   {To show this, let $\{\bar{\bx}({t_k})\}$ be a sub-sequence 
that $\lim_{k\to\infty}F(\bar{\bx}(t_k,\omega))=F^*$ (such a sub-sequence depends on the sample path $\omega$). 
Since $\omega\in \Omega_d$ and \[\lim\limits_{t\to\infty}\|\bar{\bx}(t,\omega)-\hat{z}\|\quad \text{exists}\]
for some $\hat{z}\in \zd$, we conclude that $\{\bar{\bx}(t,\omega)\}$ is a bounded sequences. Therefore, $\{\bar{\bx}(t_k,\omega)\}$ is also bounded and it has an accumulation point $z^*\in \Rd$ and hence, there is a sub-sequence $\{\bar{\bx}(t_{k_\tau},\omega)\}_{\tau\geq 0}$ of $\{\bar{\bx}(t_k,\omega)\}_{k\geq 0}$ such that 
\[\lim_{\tau\to\infty}\bar{\bx}(t_{k_\tau},\omega)=z^*.\]
As a result of continuity of $F(\cdot)$, we have 
\[\lim_{\tau\to\infty}F(\bar{\bx}(t_{k_\tau}))=F(z^*)=F^*\]
and hence, $z^*\in \mathcal{Z}$. Note that the point $z^*=z^*(\omega)$ depends on the sample path $\omega$. }

Since $\zd\subseteq \mathcal{Z}$ is dense, there is a sequence $\{q^*(s,\omega)\}_{s\geq 0}$ in $\zd$ such that $\lim_{s\to\infty}\|q^*(s,\omega)-z^*(\omega)\|=0$. Note that since $\omega\in \Omega_d$, $\lim_{t\to\infty}\|\bar{\bx}(t,\omega)-q^*(s,\omega)\|$ exists for all $s\geq 0$ and we have
    {\begin{align*}
        \lim_{t\to\infty}\|\bar{\bx}(t,\omega)-q^*(s,\omega)\|
        &=\lim_{t\to\infty}\|\bar{\bx}(t,\omega)-z^*(\omega)+z^*(\omega)-q^*(s,\omega)\|\cr
        &\leq\liminf_{t\to\infty}\|\bar{\bx}(t,\omega)-z^*(\omega)\|+\|q^*(s,\omega)-z^*(\omega)\|\cr
        &=\|q^*(s,\omega)-z^*(\omega)\|.
    \end{align*}}
    Therefore, we have 
    \begin{align}\label{eqn:limqstar}
        \lim_{s\to\infty}\lim_{t\to\infty}\|\bar{\bx}(t,\omega)-q^*(s,\omega)\|=0.
    \end{align}
    On the other hand, we have
    \begin{align*}
        \limsup_{t\to\infty}\|\bar{\bx}(t,\omega)-z^*(\omega)\|
        &=\limsup_{t\to\infty}\|\bar{\bx}(t,\omega)-q^*(s,\omega)+q^*(s,\omega)-z^*(\omega)\|\cr
        &\leq\limsup_{t\to\infty}\|\bar{\bx}(t,\omega)-q^*(s,\omega)\|+\|q^*(s,\omega)-z^*(\omega)\|\cr
        &=\left(\lim_{t\to\infty}\|\bar{\bx}(t,\omega)-q^*(s,\omega)\|\right)+\|q^*(s,\omega)-z^*(\omega)\|.
    \end{align*}
    Therefore, 
    \begin{align}\label{eqn:FINALLY}
        \limsup_{t\to\infty}\|\bar{\bx}(t,\omega)-z^*(\omega)\|
        &=\lim_{s\to\infty}\limsup_{t\to\infty }\|\bar{\bx}(t,\omega)-z^*(\omega)\|\cr 
        &\leq\lim_{s\to\infty}\lim_{t\to\infty}\|\bar{\bx}(t,\omega)-q^*(s,\omega)\|+\lim_{s\to\infty}\|q^*(s,\omega)-z^*(\omega)\|\cr 
        &=0,
    \end{align}
    where the last equality follows by combining \eqref{eqn:limqstar} and $\lim_{s\to\infty}\|q^*(s,\omega)-z^*(\omega)\|=0$. Note that \eqref{eqn:FINALLY}, implies that almost surely (i.e., for all $\omega\in\Omega_d$), we have 
    \[\lim_{t\to\infty}\bar{\bx}(t)=z^*(\omega)\]
    exists and it belongs to $\mathcal{Z}$.
    
   {Finally, according to Assumption \ref{assump:AssumpOnFunction} and \ref{assump:AssumpOnStepSize}, we have  
   \[\diam(\alpha(t)\bg(t))\leq 2{K}t^{-\beta}\max_{i\in[n]}L_i.\] Therefore, from lemma \ref{lem:ConsensusConveregnse}, we conclude that $\lim_{t\to\infty}\diam(\bx(t))=0$ almost surely, and hence  \[\lim_{t\to\infty}\|\bar{\bx}(t)-\bx_i(t)\|=0\quad \text{almost surely}.\]
   Since we almost surely have $\lim_{t\to\infty}\bar{\bx}(t)=z^*$ for a random vector $z^*$ supported in $\mathcal{Z}$,  we have $\lim_{t\to\infty}\bx_i(t)=z^*$ for all $i\in[n]$ almost surely and the proof is complete. \qed}

    % Since the left hand side is less the right hand side for every $s$, the left hand side is less than the limit of the right hand side. This implies that 
    % \begin{align*}
    %     &\limsup_{t\to\infty}\|\bar{\bx}(t,\omega)-z^*(\omega)\|\cr
    %     &\leq\lim_{s\to\infty}\left[\|q^*(s,\omega)-z^*(\omega)\|+\lim_{t\to\infty}\|\bar{\bx}(t,\omega)-q^*(s,\omega)\|\right]=0.
    % \end{align*}
    % This results in $\lim_{t\to\infty}\|\bar{\bx}(t,\omega)-z^*(\omega)\|=0$.
    %Finally, $P(\Omega)=1$ completes the proof.
    
\section{Conclusion and Future Research}\label{Conclusion}
In this work, we showed that the  averaging-based distributed optimization solving algorithm over dependent random networks converges to an optimal random point under the standard conditions on  the objective function and network formation that is conditionally $B$-connected.  To do so, we established a rate of convergence for the second moment of the autonomous averaging dynamics over such networks and used that to study the convergence of the sample-paths and second moments of the controlled variation of those dynamics. 
% Loosing almost surely to in-expectation has some implications including  robustness to link-failure, and  synthesis of fully distributed averaging-based algorithms.

Further extensions of the current work to non-convex settings, accelerated algorithms, and distributed online learning algorithms are of interest for future consideration on this topic. 
\appendix

{\it Proof of Lemma \ref{lem:DiamBehrouz}}:
% For the proof of the part (a), let 
% \[\bh\triangleq \arg\max_{j}\max_{i}\|\bx_i-\bx_j\|.\]
For the proof of part (a), let $\bx_{i}^{(k)}$ be the $k$th coordinate of $\bx_i$, and define the vector $y=[y^{(1)},\ldots,y^{(m)}]^T$ where \[y^{(k)}=\frac{1}{2}({u}^{(k)}+{U}^{(k)}),\] with 
${u}^{(k)}=\min_{i\in [n]}\bx_i^{(k)}$ and ${U}^{(k)}=\max_{i\in [n]}\bx_i^{(k)}$.
% \begin{align*}
%     (i^*(k),j^*(k))\in\arg\max_{(i,j)\in[n]\times[n]}\left|\bx_{i}^{(k)}-\bx_{j}^{(k)}\right|,
% \end{align*}
Therefore, for $\ell\in[n]$, we have
\begin{align}\label{eqn:ykmax}
    \|\bx_\ell-y\|&=\max_{k\in[m]}|\bx_{\ell}^{(k)}-y^{(k)}|\cr
    &\stackrel{}{=}\max_{k\in[m]}\left|\bx_{\ell}^{(k)}-\frac{1}{2}(u^{(k)}+U^{(k)})\right|\cr
    &\stackrel{(a)}\leq\max_{k\in[m]}\frac{1}{2}\left|u^{(k)}-U^{(k)}\right|\cr
    &=\frac{1}{2}\diam(\bx),
\end{align}
where $(a)$ follows from $u^{(k)}\leq \bx_{\ell}^{(k)}\leq U^{(k)}$.
%\todo[inline]{Prove that $\|\bx_\ell-y\|\leq \frac{1}{2}d(\bx)$}}
Also, we have
\begin{align*}
    \diam(A\bx)&=\max_{i,j\in[n]}\Bigg{\|}\sum_{\ell=1}^na_{i\ell}\bx_\ell-\sum_{\ell=1}^na_{j\ell}\bx_\ell\Bigg{\|}\cr
    &=\max_{i,j\in[n]}\Bigg{\|}\sum_{\ell=1}^na_{i\ell}(\bx_\ell-\bh)-\sum_{\ell=1}^na_{j\ell}(\bx_\ell-\bh)+\sum_{\ell=1}^n(a_{j\ell}-a_{i\ell})\bh\Bigg{\|}\cr
    &=\max_{i,j\in[n]}\Bigg{\|}\sum_{\ell=1}^n(a_{i\ell}-a_{j\ell})(\bx_\ell-\bh)\Bigg{\|}
\end{align*}
where the last equality holds as $A$ is a row-stochastic matrix  and hence, $\sum_{\ell=1}^n(a_{j\ell}-a_{i\ell})=0$. Therefore, 
\begin{align*}
    \diam(A\bx)&=\max_{i,j\in[n]}\Bigg{\|}\sum_{\ell=1}^n(a_{i\ell}-a_{j\ell})(\bx_\ell-\bh)\Bigg{\|}\cr
    &\stackrel{(a)}{\leq}\max_{i,j\in[n]}\sum_{\ell=1}^n|a_{i\ell}-a_{j\ell}|{\|}\bx_\ell-\bh{\|}\cr
    &\stackrel{(b)}{\leq}\max_{i,j\in[n]}\frac{1}{2}\diam(\bx)\sum_{\ell=1}^n|a_{i\ell}-a_{j\ell}|\cr
    &\leq\diamm(A)\diam(\bx),
\end{align*}
where $(a)$ follows from the triangle inequality, and $(b)$ follow from \eqref{eqn:ykmax}.
For the part (b), we have
\begin{align*}
    \diam(\bx+\by)&=\max_{i,j\in[n]}\|(\bx+\by)_i-(\bx+\by)_j\|\cr
    &=\max_{i,j\in[n]}\|\bx_i-\bx_j+\by_i-\by_j\|\cr
    &\stackrel{(b)}{\leq}\max_{i,j\in[n]}\left(\|\bx_i-\bx_j\|+\|\by_i-\by_j\|\right)\cr
    &\leq\max_{i,j\in[n]}\|\bx_i-\bx_j\|+\max_{i,j\in[n]}\|\by_i-\by_j\|\cr
    &=\diam(\bx)+\diam(\by),
\end{align*}
where $(b)$ follows from the triangle inequality.

For the proof of part (c), we have
\begin{align*}
    \diamm(A)&=\max_{i,j\in[n]}\sum_{\ell=1}^n \frac{1}{2}|a_{i\ell}-a_{j\ell}|\cr
    &=\max_{i,j\in[n]}\sum_{\ell=1}^n \left(\frac{1}{2}(a_{i\ell}+a_{j\ell})-\min\{a_{i\ell},a_{j\ell}\}\right)\cr
    &\stackrel{(a)}{=}\max_{i,j\in[n]}1-\sum_{\ell=1}^n \min\{a_{i\ell},a_{j\ell}\}\cr
    &=1-\min_{i,j\in[n]}\sum_{\ell=1}^n \min\{a_{i\ell},a_{j\ell}\}\cr
    &=1-\Lambda(A),
\end{align*}
where $(a)$ follows from the fact that $A$ is row-stochastic.
The proof of part (d) follows from part (c) and Lemma \ref{lem:lambdadiam}.
% , considering ${\bf b}=[B_{(1)},\ldots,B_{(n)}]^T$ where $B_{(i)}$ is the $i$th  row of a matrix $B$, and 1-norm for $\diam(\cdot)$, we have $\diamm(B)=\diam({\bf b})$. Therefore, part (a) completes the proof.

For the part (e), due to the convexity of  $\|\cdot\|$, we have 
\begin{align*}
\left\|\bx_i-\sum_{j=1}^{n}\pi_j\bx_j\right\|\leq \sum_{j=1}^{n}\pi_j\left\|\bx_i-\bx_j\right\|\leq \sum_{j=1}^{n}\pi_j\diam(\bx)=\diam(\bx).
\end{align*}
\qed %The proof of part (d) follows from $0\leq a_{ij}\leq 1$.

{\it Proof of Lemma \ref{lem:ProductwithHistory}}:
We prove by induction on $K_2$. By the assumption, the lemma is true for $K_2=K_1$. For $K_2>K_1$, we have
\begin{align*}
    \E\left[\prod_{k=K_1}^{K_2+1}Y(k)\bigg{|}\F(K_1-1)\right]
    &=\E\left[\E\left[\prod_{k=K_1}^{K_2+1}Y(k)\bigg{|}\F(K_2)\right]\bigg{|}\F(K_1-1)\right]\cr
    &=\E\left[\E\left[Y(K_2+1){|}\F(K_2)\right]\prod_{k=K_1}^{K_2}Y(k)\bigg{|}\F(K_1-1)\right]\cr
    &\leq\E\left[a(K_2+1)\prod_{k=K_1}^{K_2}Y(k)\bigg{|}\F(K_1-1)\right]\cr
    &\leq \prod_{k=K_1}^{K_2+1}a(k).
\end{align*} \qed

{\it Proof of Lemma \ref{lem:approxfunction}}:
 Consider $\hat{\tau}\geq t_0$ such that  $\hat{\theta}\triangleq\sup_{t\geq\hat{\tau}}\frac{\beta({t})}{\beta({t}+1)}\theta<1$ and let $D(t)\triangleq \sum_{s=\tau}^{t-1}\beta(s)\theta^{t-s}$. Dividing both sides by $\beta(t)>0 $, for $t> \hat{\tau}$, we have
\begin{align*}
    \frac{D(t)}{\beta(t)}&=\sum_{s=\tau}^{t-1}\frac{\beta(s)}{\beta(t)}\theta^{t-s}\cr
    &=\sum_{s=\tau}^{t-1}\prod_{\kappa=s}^{t-1}\frac{\beta(\kappa)}{\beta(\kappa+1)}\theta\cr
    &\leq\sum_{s=\tau}^{\hat{\tau}-1}\prod_{\kappa=s}^{t-1}\frac{\beta(\kappa)}{\beta(\kappa+1)}\theta+\sum_{s=\hat{\tau}}^{t-1}\hat{\theta}^{t-s}\cr
    &=\sum_{s=\tau}^{\hat{\tau}-1}\prod_{\kappa=s}^{t-1}\frac{\beta(\kappa)}{\beta(\kappa+1)}\theta+\sum_{k=1}^{t-\hat{\tau}}\hat{\theta}^{k}\cr
    &\leq \sum_{s=\tau}^{\hat{\tau}-1}\prod_{\kappa=s}^{t-1}\frac{\beta(\kappa)}{\beta(\kappa+1)}\theta+\frac{\hat{\theta}}{1-\hat{\theta}}.
\end{align*}
Let \[M_1\triangleq \sup_{t>\hat{\tau}}\sup_{\hat{\tau}\geq \tau\geq t_0}\sum_{s=\tau}^{\hat{\tau}-1}\prod_{\kappa=s}^{t-1}\frac{\beta(\kappa)}{\beta(\kappa+1)}\theta+\frac{\hat{\theta}}{1-\hat{\theta}}.\]
Note that 
\begin{align*}
    \lim_{t\to\infty}\prod_{\kappa=s}^{t-1}\frac{\beta(\kappa)}{\beta(\kappa+1)}\theta\leq\lim_{t\to\infty}\hat{\theta}^{t-\hat{\tau}}\prod_{\kappa=s}^{\hat{\tau}-1}\frac{\beta(\kappa)}{\beta(\kappa+1)}\theta=0.
\end{align*}
Therefore, $\sup_{t\geq\tau}\prod_{\kappa=s}^{t-1}\frac{\beta(\kappa)}{\beta(\kappa+1)}\theta<\infty$, and hence, $M_1<\infty$. Thus, $D(t)\leq \max\{M_1,M_2\}\beta(t)$, where $M_2\triangleq \max_{\hat{\tau}\geq t\geq \tau\geq t_0}\sum_{s=\tau}^{t-1}\frac{\beta(s)}{\beta(t)}\theta^{t-s}$, and the proof is 
complete.\qed

\bibliographystyle{abbrv}
\bibliography{reference}
\end{document}